\documentclass[a4paper, 11pt, reqno]{amsart}
\usepackage{amsmath}
\usepackage{amssymb}
\usepackage{amsthm}
\usepackage{amsxtra}
\usepackage{bm}
\usepackage[all]{xy}
\usepackage{mathrsfs}
\usepackage{mathtools}
\usepackage{stmaryrd}
\usepackage{hyperref}
\usepackage[noadjust]{cite}

\usepackage{tikz,graphicx}
\usepackage{tkz-euclide}
\usetikzlibrary{arrows,calc,graphs,quotes}
\tikzset{hidden/.style = {thick, dashed}}

\numberwithin{equation}{section}
\newtheorem{thm}{Theorem}[section]
\newtheorem{defn}[thm]{Definition}
\newtheorem{lem}[thm]{Lemma}
\newtheorem{prop}[thm]{Proposition}
\newtheorem{cor}[thm]{Corollary}

\theoremstyle{remark}
\newtheorem{rmk}[thm]{Remark}

\theoremstyle{definition}
\newtheorem{exa}[thm]{Example}

\newcommand \grp[2]{\left \langle #1 \ \middle \vert \ #2 \right \rangle}

\begin{document}

\title{Categorifications of noncrossing partitions of type B and affine type A}
\author{Yiyu Cheng}
\address{School of Mathematical Sciences, Anhui University, Hefei 230601, Anhui, PR China}
\email{cyy0327@mail.ustc.edu.cn}

\begin{abstract}
   The lattice of wide subcategories of a tube is isomorphic to the lattice of noncrossing partitions of type $B$. We prove that Kreweras complement corresponds to taking left perpendicular categories via this categorification. For a suitable choice of Coxeter element, a poset of noncrossing partition of type $\tilde A_n$ can be a lattice. We observe an inclusion from the lattice $\mathrm{NC}^B(n)$ to this lattice of type $\tilde A_n$ using weighted projective lines and prove that this inclusion preserves Kreweras complement.
\end{abstract}

\maketitle

\section{Introduction}

Noncrossing partitions are very elementary in combinatorics, and were first introduced by \cite{kreweras1972partitions} in 1972. It is related to many branches of mathematics, such as group theory, cluster theory, representation theory of algebras, and free probability. The type $B$ noncrossing partitions were due to \cite{reiner1997non} in 1997. Later, the notion of noncrossing partitions were generalized in the context of Coxeter groups by \cite{brady2001partial,brady2002k}, and independently by \cite{bessis2003dual}, see also \cite{armstrong2009generalized}.

A categorification of noncrossing partitions has been discoverd by \cite[Theorem~3.13, Theorem~3.16]{ingalls2009noncrossing}: Let $Q$ be a quiver of finite or affine type, then the poset of wide subcategories of $kQ \text{-} \mathrm{mod}$, which is ordered by inclusion, is isomorphic to the poset of non-crossing partitions of the Weyl group of $kQ$. And the theorem above was generalized to any hereditary algebra in \cite[Theorem~1.2]{hubery2016categorification}. Moreover, though this categorification, the Kreweras complement in the poset of noncrossing partitions corresponds to taking left perpendicular category in corresponding category of modules. We refer to \cite{ringel2016catalan} for an overview of Catalan combinatorics of the Hereditary Artin algebras, which illustrates the Ingalls-Thomas and Igusa-Schiffler bijections and their relations to torsion triples and stability conditions, and the Catalan combinatorics.

\cite{coelho2022functorially} works in negative cluster categories $C_{-w}(kQ)$. The authors establish a bijection between $w$-simple-minded systems in $C_{-w}(kQ)$ and positive $w$-noncrossing partitions of the corresponding Weyl group, which generalizes the Ingalls-Thomas cluster-tilting perspective to more general Calabi-Yau dimensions. And they establish a bijection between simple-minded collections in $\mathcal D^b(kQ \text{-} \mathrm{mod})$ lying in a fundamental domain of $C_{-w}(kQ)$ and $w$-simple-minded systems in $C_{-w}(kQ)$. This generalizes the same result of \cite{iyama2023positive} in the case that $Q$ is Dynkin.

\cite[Corollary 9.2(2)]{krause2021category} pointed out that the poset $\mathrm{Wid}(\mathcal T_n)$ of wide subcategories in an $n$-tube is isomorphic to the lattice $\mathrm{NC}^B(n)$ of noncrossing partitions of type $B$. And \cite[Corollary~11.3]{hanson2025para} gives a different proof by working out the representation-theory side of McCammond-Sulway \cite{mccammond2017artin}. To our knowledge, it is still unclear whether taking left perpendicular category corresponds to Kreweras complement in $\mathrm{NC}^B(n)$. In our paper, we provide an aspect of derived category for the isomorphism between $\mathrm{Wid}(\mathcal T_n)$ and $\mathrm{NC}^B(n)$. And we prove that through the isomorphism above, taking left perpendicular category corresponds to Kreweras complement. See the theorem below.

\begin{thm} \label{thm 1} (See Theorem~\ref{thm tube B})
   There is a natural isomorphism
   \[
   \phi: \mathrm{Wid}(\mathcal T_n) \cong \mathrm{NC}^B(n).
   \]
   Moreover, taking left perpendicular category corresponds to the Kreweras complement via this isomorphism, i.e. for any wide subcategory $\mathcal W$ of $\mathcal T_n$,
   \[
   \phi(\prescript{\perp}{}{\mathcal W})=K(\phi(\mathcal W)).
   \]
\end{thm}

With the help of weighted projective lines in the sense of Geigle-Lenzing, we study noncrossing partitions of type $\tilde A$ with a special Coxeter element. For a weighted projective line $\mathbb X(n)$ of weight type $(n)$, the category of coherent sheaves $\mathrm{coh} \text{-} \mathbb X(n)$ is derived equivalent with a path algebra of $\tilde{A}_n$ type. Precisely, let $Q_n$ be the quiver
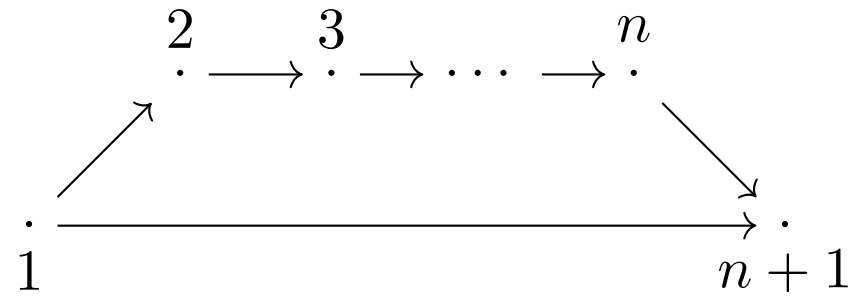
\[\tikz{
   \node (a1) at (0,0) {$\cdot$};
   \node (a2) at (1,1) {$\cdot$};
   \node (a3) at (2,1) {$\cdot$};
   \node (a4) at (3,1) {$\cdots$};
   \node (a5) at (4,1) {$\cdot$};
   \node (a6) at (5,0) {$\cdot$};

   \node (b1) at (0,-0.3) {$1$};
   \node (b2) at (1,1.3) {$2$};
   \node (b3) at (2,1.3) {$3$};
   \node (b5) at (4,1.3) {$n$};
   \node (b6) at (5,-0.3) {$n+1$};

   \graph{
           (a1)->(a2)->(a3)->(a4)->(a5)->(a6);
           (a1)->(a6);
   }
}\]
we have a derived equivalence
\[\mathcal D^b(\mathrm{coh} \text{-} \mathbb X(n)) \cong \mathcal D^b(k Q_n \text{-}\mathrm{mod}).\]
So the poset of finitely generated wide subcategories of $k Q_n \text{-} \mathrm{mod}$ is isomorphic to the poset of finitely generated wide subcategories of $\mathrm{coh} \text{-} \mathbb X(n)$, which is denoted by $\mathrm{Exc}(\mathrm{coh} \text{-} \mathbb X(n))$. Denoting by $W_{Q_n}$ the Weyl group of $k Q_n \text{-}\mathrm{mod}$ and 
\[
   c_{Q_n}=s_1s_2 \cdots s_{n+1} \in W_{Q_n}
\]
the standard Coxeter element, and $\mathrm{NC}(W_{Q_n}, c_{Q_n})$ the poset of noncrossing partitions (of $\tilde A_n$ type), from Ingalls-Thomas \cite{ingalls2009noncrossing} there is an isomorphism of posets
\[\mathrm{Exc}(\mathrm{coh} \text{-} \mathbb X(n)) \cong \mathrm{NC}(W_{Q_n}, c_{Q_n}).\]
We prove that the poset $\mathrm{Exc}(\mathrm{coh} \text{-} \mathbb X(n))$ is in fact a lattice, see Proposition~\ref{prop lattice or not}(1). And through this isomorphism, the poset of $\vec c$-invariant finitely generated wide subcategories $\mathrm{Exc}(\mathrm{coh} \text{-} \mathbb X(n))^{\vec c}$ corresponds to a sublattice 
\[B(Q_n):=\{w \in \mathrm{NC}(W_{Q_n},c_{Q_n}) \ | \ |\mathrm{Upp}(w) \cup \mathrm{Low}(w)|=\infty \}\]
of $\mathrm{NC}(W_{Q_n}, c_{Q_n})$, which is isomorphic to the lattice of noncrossing partitions of type $B$. Here, the notation $\mathrm{Upp}(w)$ means the set of upper bounds of $w$, and $\mathrm{Low}(w)$ means the set of lower bounds of $w$. Moreover, from Theorem~\ref{thm 1} we prove that the injection $B(Q_n) \to \mathrm{NC}(W_{Q_n}, c_{Q_n})$ preserves Kreweras complement. See the next theorem. 

\begin{thm} (See Theorem~\ref{thm B_n embedded in A_n tilde})
   The following statements hold.

   (1) There is a commutative diagram of lattices
   \[
   \xymatrix{\mathrm{Exc}(\mathrm{coh} \text{-} \mathbb X(n))^{\vec c} \ar@{^(->}[r] \ar[d]^{\cong} & \mathrm{Exc}(\mathrm{coh} \text{-} \mathbb X(n)) \ar[d]^{\cong} \\
   B(Q_n) \ar@{^(->}[r] & \mathrm{NC}(W_{Q_n},c_{Q_n}).
   }
   \]

   (2) There is an embedding of lattices 
   \[
     i: \mathrm{NC}^B(n) \hookrightarrow \mathrm{NC}(W_{Q_n},c_{Q_n})
   \] 
   which preserves Kreweras complement in the two lattices separately.
\end{thm}

The existence of inclusions from noncrossing partitions of type $B$ to noncrossing partitions of any affine type can be concluded by \cite[Proposition~7.6]{mccammond2017artin}. In our case, we prove that the inclusion $i$ preserves Kreweras complement.

Let us describe the structure of this paper. In Sections 2, we recall noncrossing partitions for Coxeter groups, especially for noncrossing partitions of type $A$ and type $B$. 

In Section 3, we recall some basic properties of hereditary Abelian categories, especially tubes. And we recall a categorification of noncrossing partitions due to Ingalls-Thomas \cite{ingalls2009noncrossing} and Hubery-Krause \cite{hubery2016categorification}.

In Section 4, we show a different description of Ingalls-Thomas  categorification of type $A$ which is constructed by Krause \cite{krause2021category}. It is helpful to study categorification of noncrossing partitions of type $B$. We prove that the lattice of wide subcategories of a tube is isomorphic to the lattice of noncrossing partitions of type $B$ using the language of derived category. And we prove that under this categorification, the Kreweras complement corresponds to taking left perpendicular category of the corresponding wide subcategory.

In Section 5, with the help of weighted projective line, we study noncrossing partitions of type $\tilde A_n$ associated to a specific Coxeter element. We find a sublattice of it which is isomorphic to the lattice of noncrossing partitions of type $B_n$ and it corresponds to $\vec c$-invariant finitely generated wide subcategories of $\mathrm{coh} \text{-} \mathbb X(n)$. Moreover, the sublattice induces an inclusion between lattices of noncrossing partitions which preserves Kreweras complement.

In this paper, all categories are supposed to be essentially small, i.e. the isomorphism classes of objects form a set. In the case of a weighted projective line, the base field $k$ is assumed to be algebraically closed. But in general, we do not assume a field $k$ to be algebraically closed.

\section{Noncrossing partition}

Recall that a \emph{Coxeter system} is a pair $(W,S)$, where $W$ is a group and $S$ its finite generating set, with the generating relation 
\[(st)^{m_{st}}=1,\]
where $m_{ss}=1, m_{st}=m_{ts} \geq 2$ for any $s,t \in S, s\neq t$. We say $m_{st}=\infty$ if for any $n \in \mathbb N$, $(st)^n \neq 1$. We call $|S|$ the \emph{rank} of the Coxeter system, $W$ a \emph{Coxeter group}, and $c=s_1 s_2 \cdots s_n$ a \emph{Coxeter element}, where ${s_1, \cdots, s_n}=S$. By definition, Coxeter element is not unique, and related to the generating set $S$. 

Note that a group may have different Coxeter system. For example, the dihedral group $D_{12}$ of order 12 has the following two presentations as a Coxeter group
\[
    \langle s,t \ | \ s^2=t^2=(st)^6=1 \rangle,
\]
\[
    \langle s,u,v \ | \ s^2=u^2=v^2=(su)^3=(sv)^2=(uv)^2=1 \rangle.
\]
Note that here $t=uv, v=(st)^3$.

For a Coxeter system $(W,S)$, one defines the subset of \emph{reflections} by
\[T=\{wsw^{-1} \ | \ s \in S, w \in W \}.\]
For any $w \in W$, there exists a smallest natural number $r$ such that $w=t_1 t_2 \cdots t_t$, where $t_1,\cdots,t_r \in T$. The number $r$ is called the \emph{absolute length} of $w$, denoted by $l(w)=r$. And $t_1 t_2 \cdots t_r$ is called a \emph{reduced factorization} of $w$. One defines the \emph{absolute order} on $W$ by
\[
    u \leq v \ \ \text{if} \ \ l(v) = l(u) + l(u^{-1} v).
\]
Note that $u \leq v$ if and only if there exists a reduced factorization $v=t_1 t_2 \cdots t_r$ and $p \leq r$, such that $u=t_1 t_2 \cdots t_p$. For this reason, $u$ is called a \emph{prefix} of $v$. The Coxeter group $W$ is a partially ordered set under absolute order. 

\begin{defn}
(1) Let $W$ be a Coxeter group and $c$ a Coxeter element. We call 
\[
  \mathrm{NC}(W,c):=\{w \in W \ | \ \mathrm{id} \leq w \leq c \}
\]
the poset of \emph{noncrossing partitions}.

(2) The following anti-isomorphism of the lattice $\mathrm{NC}(W,c)$ is called the \emph{Kreweras complement} 
\begin{align*}
   K: \mathrm{NC}(W,c) &\to \mathrm{NC}(W,c) \\
   w &\mapsto w^{-1}c.
\end{align*}
\end{defn}

\begin{rmk}
For a finite Coxeter group, the poset of noncrossing partitions is a lattice, but in general not.
\end{rmk}

Let $w=t_1 t_2 \cdots t_r$ be a reduced factorization, then for any subsequence $t_{i_1},\cdots,t_{i_p}$ where $i_1 < \cdots <i_p$, $t_{i_1} \cdots t_{i_p}$ is a prefix of $w$. To prove it, we only need to notice that the \emph{braid group} 
\[
   B_r=\grp{\sigma_1, \sigma_2, \cdots, \sigma_{r-1}}{
      \begin{aligned}
         \sigma_i \sigma_j &=\sigma_j \sigma_i &&\text{for } |i-j|>1, \\
         \sigma_i \sigma_{i+1} \sigma_i &= \sigma _{i+1} \sigma_i \sigma_{i+1} &&\text{for } 1 \leq i \leq r-2
      \end{aligned}
   }
\] 
has the following \emph{Hurwitz action} on the set $G^r$, where $G$ is any group
\[
    \sigma_i(g_1, \cdots, g_r):=(g_1, \cdots, g_{i-1},g_{i+1},g_{i+1}^{-1} g_i g_{i+1}, g_{i+2}, \cdots, g_r).
\]
This action first appears in \cite{hurwitz1901anzahl} for studying symmetric groups.

For Coxeter groups, the Hurwitz action on reduced factorizations of the Coxeter element is transitive. This is important for the categorification of noncrossing partitions.

\begin{prop} [{\cite[Theorem 1.4]{igusa2010exceptional}}]
   Let $(W,S)$ be a rank $n$ Coxeter system, $c$ a Coxeter element. Then $l(c)=n$ and the Hurwitz action of $B_n$ on the set of sequences $(t_1, \cdots, t_n)$ of reflections such that $t_1 \cdots t_n=c$ is transitive.
\end{prop}

\begin{rmk}
   For a Coxeter group, any subgroup generated by reflections is also a Coxeter group, see \cite[Theorem]{deodhar1989note} and \cite[Theorem 3.3]{dyer1990reflection}. So for any $w \in \mathrm{NC}(W,c)$ and reduced factorization $w=t_1 \cdots t_r$, the group $\langle t_1, \cdots, t_r \rangle$ is a Coxeter group. By the transitively of Hurwitz action, we know that the group is only determined by the element $w$.
\end{rmk}

Let us recall noncrossing partition of type $A$, which is the classical noncrossing partition \cite{kreweras1972partitions}. Let $[n]=\{1,2,\cdots,n\}$, $P=\{P_\alpha\}_{\alpha \in \Lambda}$ is a family of pairwise disjoint non-empty subsets of $[n]$. Each $P_{\alpha}$ is called a \emph{block} of $P$. We say $P$ is a \emph{partition} of $[n]$, if $\cup_{\alpha \in \Lambda} P_\alpha=[n]$. Two subsets of $[n]$ are called \emph{crossing}, if there exist $a<b<c<d$ such that $a,c$ lie in one of the subsets and $b,d$ lie in the other. 
A partition $P$ is called \emph{noncrossing}, if any two different blocks of $P$ are not crossing. The following are examples of noncrossing and crossing partitions of $[6]$.

\[
\begin{tikzpicture}
\begin{scope}[xshift=-3cm]
\draw (0,0) circle (1.5);
\node (a) at (90:1.8) {1};
\node (b) at (30:1.8) {2};
\node (c) at (-30:1.8) {3};
\node (d) at (-90:1.8) {4};
\node (e) at (-150:1.8) {5};
\node (f) at (150:1.8) {6};
\node (g) at (-90:2.5) {\{1,2,4\},\{3\},\{5,6\}};
\node (h) at (-90:3.2) {noncrossing};
\fill (90:1.5) circle (1pt);
\fill (30:1.5) circle (1pt);
\fill (-30:1.5) circle (1pt);
\fill (-90:1.5) circle (1pt);
\fill (-150:1.5) circle (1pt);
\fill (150:1.5) circle (1pt);
\filldraw[fill=gray, fill opacity=0.5] (90:1.5) to (30:1.5) to (-90:1.5) to (90:1.5);
\draw (150:1.5) to (-150:1.5);
\end{scope}
\begin{scope}[xshift=3cm]
\draw (0,0) circle (1.5);
\node (a) at (90:1.8) {1};
\node (b) at (30:1.8) {2};
\node (c) at (-30:1.8) {3};
\node (d) at (-90:1.8) {4};
\node (e) at (-150:1.8) {5};
\node (f) at (150:1.8) {6};
\node (g) at (-90:2.5) {\{1,2,4\},\{3,5,6\}};
\node (h) at (-90:3.2) {crossing};
\fill (90:1.5) circle (1pt);
\fill (30:1.5) circle (1pt);
\fill (-30:1.5) circle (1pt);
\fill (-90:1.5) circle (1pt);
\fill (-150:1.5) circle (1pt);
\fill (150:1.5) circle (1pt);
\filldraw[fill=gray, fill opacity=0.5] (90:1.5) to (30:1.5) to (-90:1.5) to (90:1.5);
\filldraw[fill=gray, fill opacity=0.2] (150:1.5) to (-150:1.5) to (-30:1.5) to (150:1.5);
\end{scope}
\end{tikzpicture}
\]

Denote by $\mathrm{NC}^A(n)$ the set of noncrossing partitions of $[n]$. It is well-known that 
\[
   |\mathrm{NC}^A(n)|=\frac{1}{n+1}\binom{2n}{n}
\]
as Catalan number, see \cite[Corollary~4.2]{kreweras1972partitions} or \cite[Theorem~4.1.2]{armstrong2009generalized}. Moreover, $\mathrm{NC}^A(n)$ has a poset structure given by refinement. Exactly, we say $\{P_\alpha \}_{\alpha \in \Lambda} \preceq \{Q_\beta \}_{\beta \in \Gamma}$ if for any $P_\alpha$, $\alpha \in \Lambda$, there exists $\beta \in \Gamma$, such that $P_\alpha \subseteq Q_\beta$. In fact, if $\{P_\alpha \}_{\alpha \in \Lambda} \preceq \{Q_\beta \}_{\beta \in \Gamma}$, then each $Q_\beta$ is the union of some $P_\alpha$. Furthermore, $\mathrm{NC}^A(n)$ is a lattice. For noncrossing partitions $P=\{P_\alpha \}_{\alpha \in \Lambda}$ and $Q=\{Q_\beta \}_{\beta \in \Gamma}$, their join is defined as follows
\[
P \wedge Q=\{P_\alpha \cap Q_\beta \ | \ \alpha \in \Lambda, \ \beta \in \Gamma, \ P_\alpha \cap Q_\beta \neq \emptyset \}.
\]
The join in the lattice of noncrossing partitions is the same with the join in the lattice of partitions. But the meet in two lattices are different.

The symmetric group $S_n$ is a Coxeter group of $A_{n-1}$ type
\[
    S_n=\langle \sigma_1, \cdots, \sigma_{n-1} \ | \ (\sigma_{i} \sigma_{i+1})^3=1, \ i=1,\cdots,n-2 \rangle,
\]
where $\sigma_1=(12)$, $\sigma_2=(23)$, $\cdots$, $\sigma_{n-1}=(n-1 \ n)$. Let $c=\sigma_1 \cdots \sigma_{n-1}=(12\cdots n)$ be the standard Coxeter element. For any $\pi \in S_n$, let $\{\pi\}$ denote the partition of $[n]$ consisting of the orbits of $\pi$. Then we have an isomorphism of lattices, see \cite[Theorem 4.1.3]{armstrong2009generalized}:
\begin{align*}
    \mathrm{NC}(S_n,c) &\xrightarrow{\sim} \mathrm{NC}^A(n) \\
    \pi &\mapsto \{ \pi \}.
\end{align*}

We recall the classical Kreweras complement for noncrossing partitions of type $A$. Let $P \in \mathrm{NC}^A(n)$, one defines its Kreweras complement $K(P)$ as follows. Let the partition $Q$ of $\{1',2',\cdots,n' \}$ be the coarsest partition such that $P \cup Q$ is a noncrossing partition of the sequence $1,1',2,2',\cdots,n,n'$. Via the map $i' \mapsto i$, $i=1,2,\cdots,n$, we obtain $K(P)$ by $Q$. Below is an example states that $K(\{1,2,4\},\{3\},\{5,6\})=\{\{1\},\{2,3\},\{4,6\},\{5\}\}$ in $\mathrm{NC}^A(6)$. 

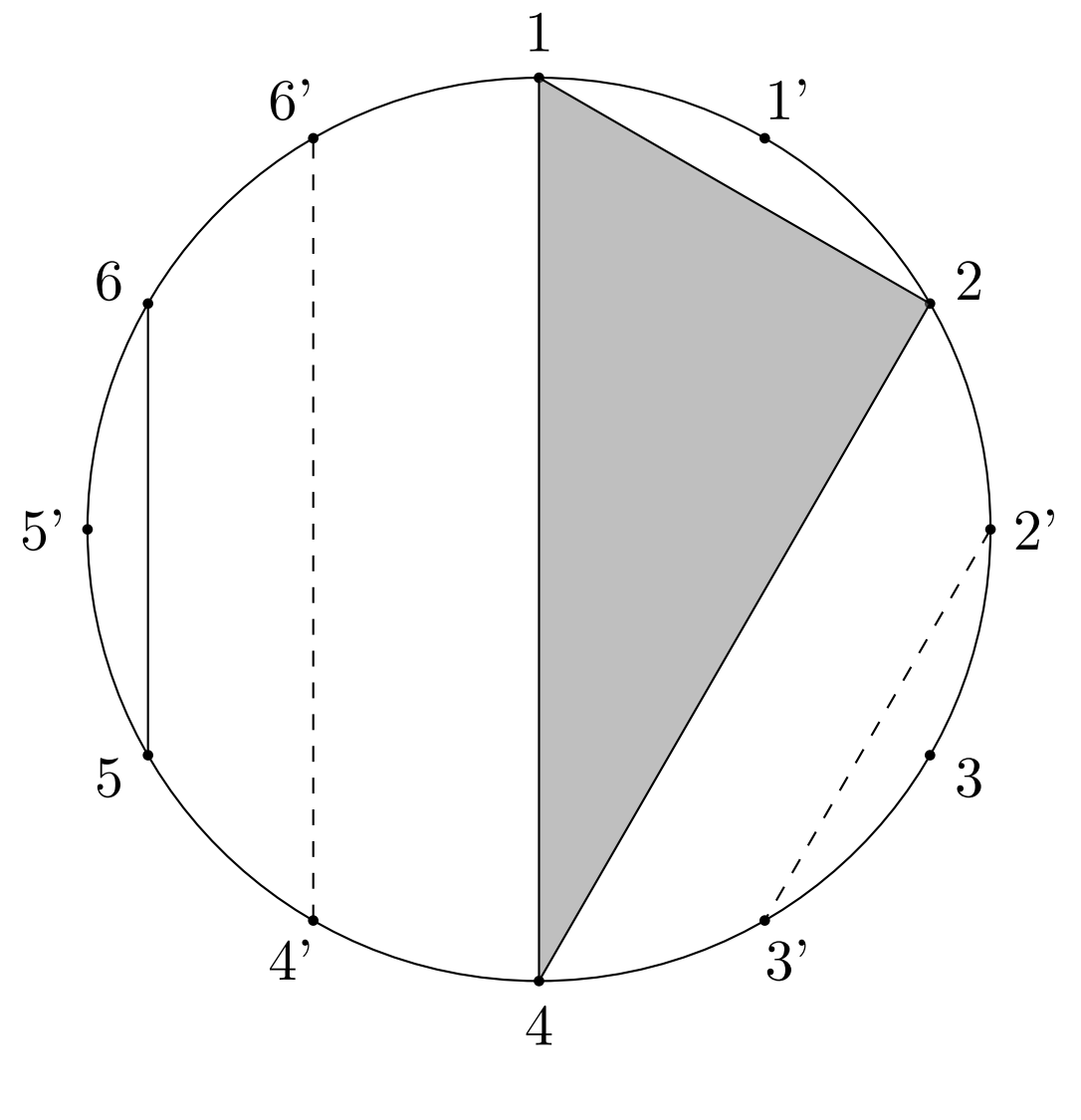
\[
\begin{tikzpicture}
\draw (0,0) circle (3);
\node (a) at (90:3.3) {1};
\node (b) at (30:3.3) {2};
\node (c) at (-30:3.3) {3};
\node (d) at (-90:3.3) {4};
\node (e) at (-150:3.3) {5};
\node (f) at (150:3.3) {6};
\node (a') at (60:3.3) {1'};
\node (b') at (0:3.3) {2'};
\node (c') at (-60:3.3) {3'};
\node (d') at (-120:3.3) {4'};
\node (e') at (180:3.3) {5'};
\node (f') at (120:3.3) {6'};

\fill (90:3) circle (1pt);
\fill (30:3) circle (1pt);
\fill (-30:3) circle (1pt);
\fill (-90:3) circle (1pt);
\fill (-150:3) circle (1pt);
\fill (150:3) circle (1pt);
\fill (60:3) circle (1pt);
\fill (0:3) circle (1pt);
\fill (-60:3) circle (1pt);
\fill (-120:3) circle (1pt);
\fill (180:3) circle (1pt);
\fill (120:3) circle (1pt);

\filldraw[fill=gray, fill opacity=0.5] (90:3) to (30:3) to (-90:3) to (90:3);
\draw (150:3) to (-150:3);

\draw[dashed] (0:3) to (-60:3);
\draw[dashed] (-120:3) to (120:3);

\end{tikzpicture}
\]

Notice that $(23)(46)=((124)(56))^{-1}(123456)$ in $S_6$. In general, \cite[Section 3]{kreweras1972partitions} points out that the classical Kreweras complement is the same with the group-theoretical Kreweras complement.

Now recall the noncrossing partition of type $B$, we refer to \cite{reiner1997non}. Consider the set
\[
    [2n]=\{ 1,2,\cdots,n,\bar 1,\bar 2,\cdots,\bar n \},
\]
where $\bar x$ is just $x+n$ (for $1 \leq x \leq n$), and assume that $\bar{\bar x}:=x$. One defines $P=(P_\alpha)_{\alpha \in \Lambda}$ as noncrossing partition of type $B$ and of order $n$ if $P$ is a noncrossing partition of $[2n]$, and for any $\alpha \in \Lambda$, $P_\alpha$ and $\overline{ P_{\alpha}}$ are both blocks of $P$. We say $P_{\alpha}$ is a \emph{zero block}, if there exists $i\in \{1, \cdots, n \}$ such that $\{i,\bar i \}\subseteq P_{\alpha}$. By the property of noncrossing, zero block is unique if it exists. We denote by $\mathrm{NC}^B(n)$ the lattice of noncrossing partitions of type $B$ of order $n$. It is well known that $|\mathrm{NC}^B(n)|=\binom{2n}{n}$, see \cite[Proposition 6]{reiner1997non}, \cite[Section 1.4]{biane2003non}. The following are two noncrossing partitions in $\mathrm{NC}^B(3)$. 
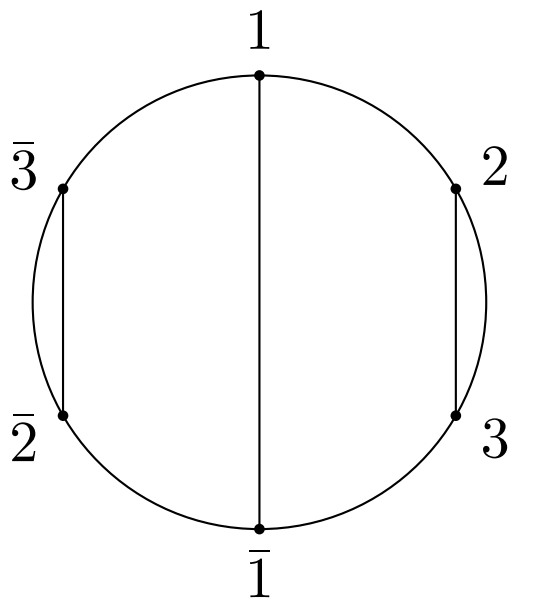
\[
\begin{tikzpicture}
\begin{scope}[xshift=-3cm]
\draw (0,0) circle (1.5);
\node (a) at (90:1.8) {1};
\node (b) at (30:1.8) {2};
\node (c) at (-30:1.8) {3};
\node (d) at (-90:1.8) {$\bar 1$};
\node (e) at (-150:1.8) {$\bar 2$};
\node (f) at (150:1.8) {$\bar 3$};
\node (g) at (-90:2.5) {$\{1,\bar 1\},\{2, 3\},\{\bar 2, \bar 3\}$};
\node (h) at (-90:3.2) {has zero block $\{1, \bar 1 \}$};
\fill (90:1.5) circle (1pt);
\fill (30:1.5) circle (1pt);
\fill (-30:1.5) circle (1pt);
\fill (-90:1.5) circle (1pt);
\fill (-150:1.5) circle (1pt);
\fill (150:1.5) circle (1pt);
\draw (90:1.5) to (-90:1.5);
\draw (30:1.5) to (-30:1.5);
\draw (150:1.5) to (-150:1.5);
\end{scope}
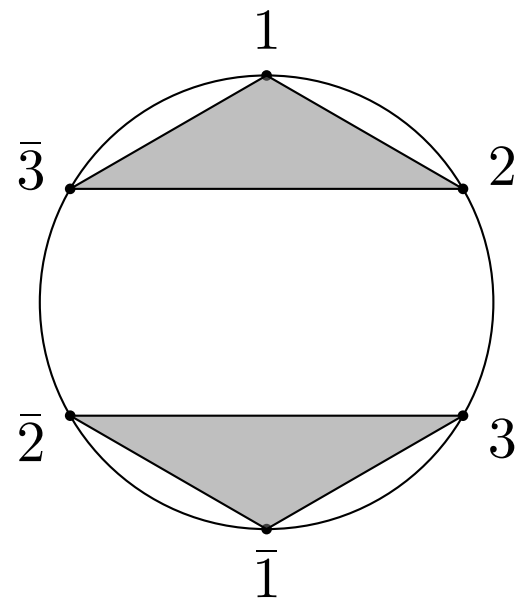
\begin{scope}[xshift=3cm]
\draw (0,0) circle (1.5);
\node (a) at (90:1.8) {1};
\node (b) at (30:1.8) {2};
\node (c) at (-30:1.8) {3};
\node (d) at (-90:1.8) {$\bar 1$};
\node (e) at (-150:1.8) {$\bar 2$};
\node (f) at (150:1.8) {$\bar 3$};
\node (g) at (-90:2.5) {$\{1,2,\bar 3\},\{3,\bar 1, \bar 2\}$};
\node (h) at (-90:3.2) {has no zero block};
\fill (90:1.5) circle (1pt);
\fill (30:1.5) circle (1pt);
\fill (-30:1.5) circle (1pt);
\fill (-90:1.5) circle (1pt);
\fill (-150:1.5) circle (1pt);
\fill (150:1.5) circle (1pt);
\filldraw[fill=gray, fill opacity=0.5] (150:1.5) to (90:1.5) to (30:1.5) to (150:1.5);
\filldraw[fill=gray, fill opacity=0.5] (-150:1.5) to (-90:1.5) to (-30:1.5) to (-150:1.5);
\end{scope}
\end{tikzpicture}
\]

Now we explain the relationship between noncrossing partitions of type $B$ and Coxeter group $W(B_n)$ of type $B$, we refer to \cite[Section~4.5]{armstrong2009generalized}. Recall the Coxeter diagram of type $B$
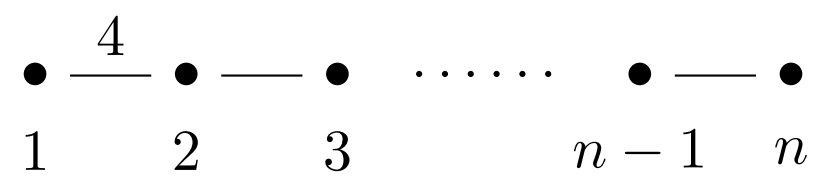
\[
    \tikz{
            \node (a_1) at (0,0) {$\bullet$};
            \node (a_2) at (1,0) {$\bullet$};
            \node (a_3) at (2,0) {$\bullet$};
            \node (a_4) at (3,0) {$\cdots \cdots$};
            \node (a_5) at (4,0) {$\bullet$};
            \node (a_6) at (5,0) {$\bullet$};
            \node (b_1) at (0,-0.5) {1};
            \node (b_2) at (1,-0.5) {2};
            \node (b_3) at (2,-0.5) {3};
            \node (b_5) at (4,-0.5) {$n-1$};
            \node (b_6) at (5,-0.5) {$n$};

            \graph[nodes={circle, fill}]{
            (a_1) --["4"] (a_2) --(a_3);
            (a_5) --(a_6);};
         }
\]
Corresponding to the points in the Coxeter graph above, we denote by $s_1, s_2, \cdots, s_n$ the generators of $W(B_n)$, and let $c=s_1 s_2 \cdots s_n$ be the standard Coxeter element. There is a natural inclusion of groups
\begin{align*}
   i_B: W(B_n) &\to S_{2n}, \\
   s_1 &\mapsto (1 \bar 1), \\
   s_2 &\mapsto (1 2)(\bar 1 \bar 2), \\ 
   \cdots \\
   s_n &\mapsto (n-1 \ n)(\overline{n-1} \ \overline n). \\
\end{align*}
The image of this inclusion is consisting of all maps which  commute with the \emph{antipodal map} $a=(1 \bar 1)(2 \bar 2) \cdots (n \bar n)$, i.e. $\mathrm{im}(i_B)=C_{S_{2n}}(a)$, the centralizer of $a$, see \cite[Section~4.5]{armstrong2009generalized}. Note that $i_B$ preserves the standard Coxeter element
\begin{align*}
   i_B(c) &=i_B(s_1) i_B(s_2) \cdots i_B(s_n) \\
    &=(1 \bar 1)(12 \cdots n)(\bar 1 \bar 2 \cdots \bar n) \\
    &=(12 \cdots n \bar 1 \bar 2 \cdots \bar n),
\end{align*}
and $a=(1 \bar 1)(2 \bar 2) \cdots (n \bar n)=i_B(c^n)$. Through the map $i_B$, we get the following isomorphism of lattices, see \cite[Theorem 3.2]{biane2003non}
\begin{align*}
   \mathrm{NC}(W(B_n),c) &\xrightarrow{\sim} \mathrm{NC}(S_{2n},i_B(c)) \cap C_{S_{2n}}(a),\\
   w &\mapsto i_B(w).
\end{align*}
Similar with the case of type $A$, there are isomorphisms of lattices, see \cite[Theorem 4.5.3]{armstrong2009generalized}
\begin{alignat*}{5}
    \mathrm{NC}(S_{2n},i_B(c)) \cap C_{S_{2n}}(a) &\xrightarrow{\sim} \mathrm{NC}^B(n) , \\
    \pi &\mapsto \{\pi\} .
\end{alignat*}
The composition of the two isomorphisms above gives the isomorphism of lattices $\mathrm{NC}(W(B_n),c) \xrightarrow{\sim} \mathrm{NC}^B(n)$, which preserves the standard Coxeter element.

Let $P \in \mathrm{NC}^B(n)$, the Kreweras complement of $P$ is just $K(P)$ in $\mathrm{NC}^A(2n)$. It is same with the group-theoretical Kreweras complement since the map $i_B$ preserves the standard Coxeter element. 

\section{Wide subcategories of a hereditary Abelian category}

By Ingalls-Thomas \cite{ingalls2009noncrossing} and Hubery-Krause \cite{hubery2016categorification}, the poset of noncrossing partitions of a Coxeter group is isomorphic to the poset of wide subcategories of the category of finite generated modules over a hereditary algebra. We recall this categorification  of noncrossing partitions in the framework of hereditary abelian categories.

Let $\mathcal A$ be an abelian category, and $\mathcal W$ be a full subcategory of $\mathcal A$. We call $\mathcal W$ a \emph{wide} subcategory if it is closed under taking kernels, cokernels, and extensions. Note that a wide subcategory is closed under direct summands.

An abelian category $\mathcal H$ is called \emph{hereditary}, if $\mathrm{Ext}_{\mathcal H}^2(X,Y)=0$ for any $X,Y \in \mathcal H$. In this paper, we mainly consider the following 3 kinds of hereditary abelian categories: the category of modules over a hereditary algebra, the category of coherent sheaves over a weighted projective line, tubes. Note that the first two categories admit a tilting object, but the third not. And all of them are \emph{Ext-finite}, i.e. for any $X,Y \in \mathcal H$ and $i \geq 0$, we have $\mathrm{dim}_k \mathrm{Ext}_{\mathcal H}^i(X,Y)< \infty$.

For a hereditary abelian category $\mathcal H$, we denote by $\mathcal D^b(\mathcal H)$ the bounded derived category of $\mathcal H$. A triangulated subcategory $\mathcal D$ of $\mathcal D^b(\mathcal H)$ is called \emph{thick} if it is closed under taking direct summands. By the following proposition, there is a one-to-one correspondence between wide subcategories of a hereditary abelian category $\mathcal H$ and thick subcategories of $\mathcal D^b(\mathcal H)$.

\begin{prop} [{\cite[Theorem~5.1]{bruning2007thick}}] \label{bruning}
   Let $\mathcal H$ be a hereditary abelian category. Then the following statements hold.\\
   (1) For any $X \in \mathcal D^b(\mathcal H)$, we have a non-canonical isomorphism
   \[
   X \cong \bigoplus_{n \in \mathbb Z}H^nX[-n].
   \] \\
   (2) We have a one-to-one correspondence between wide subcategories of $\mathcal H$ and thick subcategories of $\mathcal D^b(\mathcal H)$, by assignments
   \[
   \mathcal W \mapsto \{X \in \mathcal D^b(\mathcal H) \; | \; H^n(X) \in \mathcal W, \ \forall n \in \mathbb Z \} \ and \  \mathcal T \mapsto \{ H^0(X) | X \in \mathcal T \}.
   \]
   \end{prop}

Let $\mathcal H$ be a hereditary abelian category, and $\mathcal X \subseteq \mathcal H$ be a class of objects. We denote by $\langle \mathcal X \rangle$ the smallest wide subcategory of $\mathcal H$ containing $\mathcal X$, which is called the wide subcategory \emph{generated} by $\mathcal X$. For a triangulated category $\mathcal T$ and $\mathcal X \subseteq \mathcal T$, we still use the notation $\langle \mathcal X \rangle$ to present the smallest thick subcategory of $\mathcal T$ containing $\mathcal X$. The right and left \emph{perpendicular categories} of $\mathcal X$ are defined to be the full subcategories
\[
\mathcal X^{\perp} :=\{A \in \mathcal H \ | \ \mathrm{Hom}_{\mathcal H}(X,A)=0=\mathrm{Ext}_{\mathcal H}^1(X,A), \ \forall X \in \mathcal X \}
\]
and
\[
\prescript{\perp}{}{\mathcal X} :=\{A \in \mathcal H \ | \ \mathrm{Hom}_{\mathcal H}(A,X)=0=\mathrm{Ext}_{\mathcal H}^1(A,X), \ \forall X \in \mathcal X \},
\]
respectively. It is well known that both $\mathcal X^{\perp}$ and $\prescript{\perp}{}{\mathcal X}$ are wide subcategories; see \cite[Proposition~1.1]{Werner1991Perpendicular}. And then we have $\mathcal X^{\perp}=\langle \mathcal X \rangle^{\perp}$ and $\prescript{\perp}{}{\mathcal X}=\prescript{\perp}{}{\langle \mathcal X \rangle}$.

For a triangulated category $\mathcal T$ and a subcategory $\mathcal X$, one defines perpendicular categories as follows.
\[
\mathcal X^{\perp} :=\{A \in \mathcal T \ | \ \mathrm{Hom}_{\mathcal H}(X,A[i])=0, \ \forall X \in \mathcal X, i\in \mathbb Z \}
\]
and
\[
\prescript{\perp}{}{\mathcal X} :=\{A \in \mathcal T \ | \ \mathrm{Hom}_{\mathcal H}(A,X[i])=0, \ \forall X \in \mathcal X, i \in \mathbb Z \}.
\]
It is clear that for a hereditary Abelian category $\mathcal H$, the correspondence between wide subcategories of $\mathcal H$ and thick subcategories of $\mathcal D^b(\mathcal H)$ in Proposition~\ref{bruning} preserves taking perpendicular categories.

\begin{defn} \label{exceptional}
Let $\mathcal A$ be an abelian $k$-category and let $E \in \mathcal A$. 

(1) We call $E$ a \emph{rigid} object if $\mathrm{Ext}_{\mathcal A}^p(E,E)=0$ for any $p \geq 1$. 

(2) We call $E$ an \emph{exceptional} object if $E$ is rigid and $\mathrm{End}_{\mathcal A}(E)$ is a division algebra. 

(3) A sequence of objects $(E_1,E_2, \cdots ,E_n)$ is called an \emph{exceptional sequence} if each $E_i$ is exceptional and $\mathrm{Ext}_{\mathcal A}^p(E_i,E_j)=0$ for any $i>j$ and $p \geq 0$.

(4) An exceptional sequence is \emph{full} provided that its objects generate $\mathcal A$.
\end{defn}

For a triangulated category $\mathcal T$, an \emph{exceptional object} $E$ is defined by 
\[
\mathrm{Hom}_{\mathcal T}(E,E[i])=\begin{cases}
   \text{a division ring}, &\text{if } i=0; \\
   0, &\text{if } i \neq 0.
\end{cases}
\]
A sequence of exceptional objects $E_1,\cdots,E_n$ is called an \emph{exceptional sequence} if $\mathrm{Hom}_{\mathcal T}(E_i,E_j[l])=0$ for any $i>j$ and $l \neq 0$.

For a hereditary abelian category $\mathcal H$, we denote by $\mathrm{Wid}(\mathcal H)$ the set of wide subcategories of $\mathcal H$, and $\mathrm{Exc}(\mathcal H)$ the set of wide subcategories which are generated by an exceptional sequence. In \cite{ingalls2009noncrossing}, a wide subcategory generated by an exceptional sequence is called \emph{finitely generated}. It is well known that in an Ext-finite hereditary Abelian $k$-category, a wide subcategory $\mathcal X$ is finitely generated if and only if there exists a rigid object $X$ such that $\mathcal X=\langle X \rangle$, see \cite[Proposition~2.7]{cheng2024wide}. We denote by 
\[\mathrm{NExc}(\mathcal H):=\mathrm{Wid}(\mathcal H) \setminus \mathrm{Exc}(\mathcal H).\] 
We view the zero category as the subcategory generated by the zero exceptional sequence. In \cite{ingalls2009noncrossing}, a wide subcategory generated by an exceptional sequence is called finitely generated.

A pair $(\mathcal X,\mathcal Y)$ of full subcategories is called \emph{Ext-orthogonal} if $\mathcal X^{\perp}=\mathcal Y$ and $\mathcal X=\prescript{\perp}{}{\mathcal Y}$. The following perpendicular calculus is useful. 

\begin{prop}[{\cite[Proposition 3.3]{cheng2024wide}}] \label{perp_cal}
   Let $\mathcal H$ be an Ext-finite hereditary abelian category and let $X \in \mathcal H$ be a rigid object. Then the following statements hold. 
   
   (1) Any $M \in \mathcal H$ fits into a 5-term exact sequence
   \[
   0 \to Y^1 \to X_0 \to M \to Y^0 \to X_1 \to 0
   \]
   with $X_0,X_1 \in \langle X \rangle$ and $Y^0,Y^1 \in X^{\perp}$. 
   
   (2) The pair $(\langle X \rangle,X^{\perp})$ is Ext-orthogonal.
\end{prop}

We might compare the condition of Proposition~\ref{perp_cal} and  \cite[Proposition~A.1]{hubery2016categorification} via the next proposition.

\begin{prop} [{\cite[Proposition 3.6]{cheng2024wide}}] \label{proj_gen}
   Let $\mathcal H$ be a Hom-finite hereditary abelian category and $P$ a projective object of $\mathcal H$. Then $\mathcal H=\langle P \rangle$ if and only if $P$ is a generator of $\mathcal H$, i.e. for any $X \in \mathcal H$, there exists $n \in \mathbb N$ and an epimorphism $P^n \twoheadrightarrow X$.
\end{prop}

An object $T$ in a hereditary Abelian category $\mathcal H$ is called a \emph{tilting object} if $T$ is rigid and $\mathcal H=\langle T \rangle$. The following Proposition to determine a perpendicular category is equivalent to a module category is known to experts.

\begin{prop} \label{prop perp cat is module cat}
   Let $\mathcal H$ be a Ext-finite hereditary Abelian $k$-category, and let $X \in \mathcal H$ be a rigid object. If there exists a tilting object $T$ such that $T$ can be embedded in some $X^n$, then $X^{\perp}$ has a projective generator, $\prescript{\perp}{}{X}$ has an injective cogenerator.
\end{prop}

\begin{proof}
   Assume that $T \hookrightarrow X^n$. Since $\mathcal H$ is hereditary, we have $\mathrm{Ext}_{\mathcal H}^1(T,X)$=0. By Bongartz complement $\mathcal H$ has a tilting object $X \oplus E$. Precisely, $E$ is the middle term of the following right universal extension
   \[0 \to T \to E \to X' \to 0,\]
   where $X' \in \mathrm{add}(X)$, the subcategory of direct summands of direct sums of $X$. By the word right universal extension one mean the connecting map 
   \[\delta: \mathrm{Hom}_{\mathcal H}(X, X') \to \mathrm{Ext}_{\mathcal H}^1(X,T) \] 
   is surjective. Since $X$ is rigid, $T \hookrightarrow X^n$ can be extended to $E \to X^n$. And then $E$ can be embedded in $X' \oplus X^n$ by a similar argument with the Horseshoe lemma. That is to say $E$ can be cogenerated by $X$. It follows that $\mathrm{Hom}_{\mathcal H}(\prescript{\perp}{}{X},E)=0$. So by a dual version of Proposition~\ref{perp_cal}, $E$ fits into a five-term exact sequence and it degenerates to a short exact sequence
   \[0 \to E \xrightarrow{f} X'' \to I \to 0,\]
   where $X'' \in \langle X \rangle$, $I \in \prescript{\perp}{}{X}$. We will show that $I$ is an injective cogenerator of  $\prescript{\perp}{}{X}$. For any $M \in \prescript{\perp}{}{X}$, applying $\mathrm{Hom}_{\mathcal H}(M,-)$ on $\epsilon$ we obtain $\mathrm{Ext}_{\mathcal H}^1(M,I)=0$. So $I$ is injective in $\prescript{\perp}{}{X}$. Obviously $\langle I,X \rangle=\langle E,X \rangle=\mathcal H$, by Proposition~\ref{perp_cal}, $\langle I \rangle=\prescript{\perp}{}{X}$. By a dual version of Proposition~\ref{proj_gen}, $I$ is an injective cogenerator of $\prescript{\perp}{}{X}$.

   Put $E$ into the five-term exact sequence given by rigid object $X$
   \[
      0 \to Y^1 \to X_0 \to E \to P \to X_1 \to 0,
   \]
   where $X_0,X_1 \in \langle X \rangle$, $P,Y^1 \in X^{\perp}$. We will prove that $P$ is a projective generator of $X^{\perp}$. Denote by $N$ the image of $E \to P$, then we have a short exact sequence
   \[
      \eta: 0 \to N \to P \to X_1 \to 0.
   \]
   Since $\mathcal H$ is hereditary, so for any $Y \in X^{\perp}$, $E$ can be cogenerated by $X$ implies that $\mathrm{Ext}_{\mathcal H}^1(E,Y)=0$. Since $N$ is a quotient object of $E$, we have $\mathrm{Ext}_{\mathcal H}^1(N,Y)=0$. Applying $\mathrm{Hom}_{\mathcal H}(-,Y)$ to $\eta$, we obtain $\mathrm{Ext}_{\mathcal H}^1(P,Y)=0$. So $P$ is projective in $X^{\perp}$, specifically, $P$ is rigid. Obviously $\langle P,X \rangle=\langle E,X \rangle=\mathcal H$. So by Proposition~\ref{perp_cal}, $\langle P \rangle=X^{\perp}$. And then by Proposition~\ref{proj_gen}, $P$ is a projective generator of $X^{\perp}$.
\end{proof}

Let $A$ be a finite dimensional hereditary algebra over a field $k$. Recall that the \emph{Cartan form} on $K_0(A \text{-mod})$ is a bilinear form
\begin{align*}
   I_{A}:K_0(A \text{-mod}) \times K_0(A \text{-mod}) &\to \mathbb Z, \\
   ([X],[Y]) &\mapsto \langle [X],[Y] \rangle + \langle [Y],[X] \rangle,
\end{align*}
where $\langle -,- \rangle$ is the \emph{Euler form} 
\[
   \langle [X],[Y] \rangle=\mathrm{dimHom}_{A}(X,Y)-\mathrm{dimExt}_{A}^1 (X,Y).
\]
Let $(E_1,E_2, \cdots, E_n)$ be a full exceptional sequence of $A\text{-mod}$, then the \emph{Weyl group} $W_A$ of $A \text{-mod}$ is defined to be the subgroup of $\mathrm{Aut}(K_0(A \text{-mod}), I_{A})$ generated by reflections
\[
   r_{[E_i]}(v):=v-I_{\mathcal H}(v,[E_i])[E_i], \ \ v\in K_0(\mathcal H), \ \ i=1,2,\cdots,n.
\]
Let $c_A=r_{[E_1]} r_{[E_2]} \cdots r_{[E_n]}$ called the \emph{Coxeter element} of the Weyl group $W_A$.

Note that simple modules of $A$ can be reordered to be a full exact sequence of $A \text{-mod}$. And there is a transitive action of braid groups on the set of full exceptional sequences \cite[Theorem]{crawley1993exceptional}. So $W_A$ is a Coxeter group and $c_A$ is a Coxeter element in the sense of Coxeter system by \cite[Proposition 3.13]{kac1990infinite}. 

The following theorem states a categorification of noncrossing partitions. Note that the map cox is well-defined because of the Hurwitz transitivity in Coxeter groups.

\begin{thm} [{\cite[Theorem 3.13, Theorem 3.16]{ingalls2009noncrossing}}{\cite[Theorem~1.2]{hubery2016categorification}}] \label{hubery_krause}

   Let $A$ be a finite dimensional hereditary algebra with Weyl group $W$ and $c_A$ the Coxeter element of $W_A$. Then
   \begin{align*}
      \mathrm{cox}: \mathrm{Exc}(A \text{-} \mathrm{mod}) &\to \mathrm{NC}(W_A,c_A) \\
      \langle E_1,\cdots, E_r \rangle &\mapsto r_{[E_1]} \cdots r_{[E_r]}
   \end{align*}
   is an isomorphism of posets.
\end{thm}

Moreover, the Kreweras complement in $\mathrm{NC}(W_A,c_A)$ corresponds to taking left perpendicular category in $A \text{-mod}$.

\begin{prop} [{\cite[Proposition 6.1]{hubery2016categorification}}] \label{Kreweras in cat} 

   Let $\mathcal W \in \mathrm{Exc}(A \text{-} \mathrm{mod})$, then 
   \[
     \mathrm{cox}(\prescript{\perp}{}{\mathcal W})=K(\mathrm{cox}(\mathcal W)).
   \]
\end{prop}

Now we recall tubes and their wide subcategories. 

Let $\mathcal U$ be a \emph{length abelian category}, that is, an abelian category whose objects have finite length. Let $E$ be an object, we denote by $l(E)$ the length of $E$. An object $U$ is called \emph{uniserial} if it has a unique composition series. If all indecomposable objects in $\mathcal U$ are uniserial, then we call $\mathcal U$ a \emph{uniserial category}.

Let $\mathcal U$ be a uniserial category and $S$ a simple object in $\mathcal U$. Up to isomorphism, we denote by $S_{[n]}$ the unique indecomposable object with socle $S$ and length $n$ (such $S_{[n]}$ may not exist). Then $S_{[n]}$ has the unique chain of all subobjects:
\[
0 \subseteq S \subseteq S_{[2]} \subseteq \cdots \subseteq S_{[n-1]} \subseteq S_{[n]}.
\]
And each indecomposable object in $\mathcal U$ is of the form $S_{[n]}$; see \cite[Subsection~1.7]{chen2009introduction}. Dually, we denote by $S^{[n]}$ the unique indecomposable object with top $S$ and length $n$. Note that the notations $S_{[n]}$ and $S^{[n]}$ are not standard.

For an object $E$ in a length category, we denote by $\mathrm{top}(E)$ and $\mathrm{soc}(E)$ the top and the socle of $E$, respectively. 

\begin{lem} \label{uni_hom}
   Let $\mathcal U$ be a uniserial category and $E$, $F$ indecomposable objects. Then $\mathrm{Hom}_{\mathcal U}(E,F) \neq 0$ if and only if $\mathrm{top}(E)$ is a composition factor of $F$ and $\mathrm{soc}(F)$ is a composition factor of $E$.
\end{lem}
   
\begin{proof}
   All objects in this proof are assumed to be indecomposable. 
   
   If $\mathrm{Hom}_{\mathcal U}(E,F) \neq 0$, we take a nonzero morphism $f: E \to F$. We have $\mathrm{top}(\mathrm{im} f) \cong \mathrm{top}(E)$ and $\mathrm{soc}(\mathrm{im} f) \cong \mathrm{soc}(F)$. So $\mathrm{top}(E)$ is a composition factor of $F$ and $\mathrm{soc}(F)$ is a composition factor of $E$. 
   
   Conversely, since $\mathrm{soc}(F)$ is a composition factor of $E$, there exists a quotient object $Q$ of $E$ whose socle is isomorphic to $\mathrm{soc}(F)$. So there exists a shortest object $I$ in $\mathcal U$ such that $\mathrm{top}(I)\cong \mathrm{top}(E)$ and $\mathrm{soc}(I)\cong \mathrm{soc}(F)$. So $I$ is a quotient of $Q$ since $\mathrm{top}(Q) \cong \mathrm{top}(I)$ and $l(Q) \geq l(I)$. Then $I$ is a quotient of $E$. 
   
   Since $\mathrm{top}(E)$ is a composition factor of $F$, there exists a subobject $Q'$ of $F$ with $\mathrm{top}(Q') \cong \mathrm{top}(E)$. By the minimality of $l(I)$, we have $l(I) \leq l(Q') \leq l(F)$. Then $\mathrm{soc}(I) \cong \mathrm{soc}(F)$ and $l(I) \leq l(F)$ implies that $I$ is a subobject of $F$.

   Now we can get a composition of an epimorphism and a monomorphism $E \twoheadrightarrow I \hookrightarrow F$, which implies that $\mathrm{Hom}_{\mathcal U}(E,F) \neq 0$.
\end{proof}

Let $\mathcal H$ be a Hom-finite abelian $k$ category. We say that $\mathcal H$ satisfies \emph{Serre duality} if there exist an equivalence $\tau: \mathcal H \to \mathcal H$ and natural isomorphisms
\[
\mathrm{Ext}_{\mathcal H}^1(X,Y) \xrightarrow{\sim} \mathrm{D \ Hom}_{\mathcal H}(Y, \tau X)
\]
for all objects $X,Y$ from $\mathcal H$. Here, $\mathrm{D}=\mathrm{Hom}_k(-,k)$ denotes the $k$-duality. It follows from \cite[Proposition~1.4]{lenzing2006hereditary} that $\mathcal H$ is hereditary without nonzero projectives or injectives. It is well known that a Hom-finite length category with Serre duality is uniserial. We refer to \cite[Theorem 1.7]{lenzing2006hereditary} and \cite[Proposition 8.3]{gabriel1973indecomposable}.

Recall that an additive category $\mathcal A$ is \emph{connected}, if any decomposition $\mathcal A=\mathcal A_1 \coprod \mathcal A_2$ into additive categories implies that $\mathcal A_1=0$ or $\mathcal A_2=0$. 

In what follows, we denote by $\mathcal T_n$ a connected uniserial hereditary abelian category with Serre duality which has exactly $n$ simples up to isomorphism, which is also called an \emph{$n$-tube}. 

\begin{rmk}
   { \rm Let $k$ be a field and $\Gamma _n$ the quiver of $\tilde{\mathbb A}_n$ type with cyclic orientation. Then the $k$-category $\mathcal U_n$ is equivalent to $\mathrm{rep}_0(\Gamma_n, k)$, the category of nilpotent representations of $\Gamma_n$. }
\end{rmk}

Note that full exceptional sequences do not exist in tubes. Moreover, taking left and right perpendicular categories give a duality between $\mathrm{Exc}(\mathcal T_n)$ and $\mathrm{NExc}(\mathcal T_n)$.

\begin{prop} [{\cite[Theorem 2.3.25]{dichev2009thick}, \cite[Proposition 4.6]{cheng2024wide}}] \label{tube perpendicular}
   Let $\mathcal W \in \mathcal T_n$. Then
   \[\prescript{\perp}{}{(W^{\perp})}=\mathcal W=(\prescript{\perp}{}{W})^{\perp}.\]
   Furthermore,
   \begin{align*}
      \mathrm{Exc}(\mathcal T_n)^{\perp}&=\mathrm{NExc}(\mathcal T_n)=\prescript{\perp}{}{\mathrm{Exc}(\mathcal T_n)} \\
      \mathrm{NExc}(\mathcal T_n)^{\perp}&=\mathrm{Exc}(\mathcal T_n)=\prescript{\perp}{}{\mathrm{NExc}(\mathcal T_n)}.
   \end{align*}
\end{prop}

\section{Tube and noncrossing partition of type B}

Let $\mathcal A$ be an Abelian category, an object $X$ is called a \emph{brick}, if $\mathrm{End}_{\mathcal A}(X)$ is a division ring. Two bricks $X$ and $Y$ are called \emph{orthogonal}, if $\mathrm{Hom}_{\mathcal A}(X,Y)=0=\mathrm{Hom}_{\mathcal A}(Y,X)$. 

Let $\mathcal X$ be a set (or a class) of orthogonal bricks, i.e. each object in $\mathcal X$ is a brick and they are pairwise orthogonal. We denote by $\mathrm{Orth}(\mathcal A)$ the set of orthogonal sets of bricks of $\mathcal A$. An \emph{$\mathcal X$-filtration} of an object $A$ is defined to be a chain of subobjects
\[
      0=A_0 \subseteq A_1 \subseteq A_2 \subseteq \cdots \subseteq A_n=A,
\]
such that $A_i/A_{i-1} \in \mathcal X$, $1 \leq i \leq n$. We denote by $\mathcal U(\mathcal X)$ the full subcategory of $\mathcal A$ formed by objects which have an $\mathcal X$-filtration. 

The following theorem is the famous ``Ringel's simplification".

\begin{thm} [{\cite[Theorem~1.2]{ringel1976representations} }]\label{thm ringel simplification}
   Let $\mathcal A$ be an Abelian category, $\mathcal X$ a set of orthogonal bricks. Then $\mathcal U(\mathcal X)$ is a wide subcategory of $\mathcal A$, i.e. $\langle \mathcal X \rangle=\mathcal U(\mathcal X)$. And $\mathcal X$ is the set of all simple objects in $\mathcal U(\mathcal X)$.
\end{thm}

\begin{cor} \label{cor brick wide cat}
   Let $\mathcal A$ be a length category, $\mathcal W$ a wide subcategory of $\mathcal A$, $\mathcal X$ a set of orthogonal bricks. Denote by $\mathcal S(\mathcal W)$ a complete set of representatives of the isomorphic classes of simple objects in $\mathcal W$. Then signs $\mathcal X \mapsto \mathcal U(\mathcal X)$ and $\mathcal W \mapsto \mathcal S(\mathcal W)$ gives a mutually inverse bijections between $\mathrm{Orth}(\mathcal A)$ and $\mathrm{Wid}(\mathcal A)$.
\end{cor}
\begin{proof}
   Firstly, any wide subcategory $\mathcal W$ is a length category, since the ascending and descending chains of subobjects in $\mathcal W$ are stable. By the definition of simple object, $\mathcal S(\mathcal W)$ is a set of orthogonal bricks of $\mathcal A$.

   Since $\mathcal W$ is a length category, each object has a decomposition series in $\mathcal W$, so $\mathcal U(\mathcal S(\mathcal W))=\mathcal W$. On the other hand, by Theorem \ref{thm ringel simplification}, $\mathcal X$ is just the set of simple objects of $\mathcal U(\mathcal X)$, that is $\mathcal S(\mathcal U(\mathcal X))=\mathcal X$.
\end{proof}

\begin{rmk}
   {\rm $\mathrm{Orth}(\mathcal A)$ is a partially ordered set given by $\mathrm{Wid}(\mathcal A)$ and the bijection above.}
\end{rmk}

Based on the above corollary, \cite[Lemma 6.2]{krause2021category} gives a new description of Ingalls-Thomas correspondence in Theorem \ref{hubery_krause} for algebras of type $A$. We now begin to illustrate it. The quiver $\vec{\mathbb A}_{n}$ is as follows:
\[
    \tikz {
        \node (a_1) at (0,0) {$\bullet$};
        \node (a_2) at (1.5,0) {$\bullet$};
        \node (a_3) at (3.25,0) {$\cdots$};
        \node (a_4) at (5.1,0) {$\bullet$};
        \node (b_1) at (0,-0.5) {$1$};
        \node (b_2) at (1.5,-0.5) {$2$};
        \node (b_4) at (5.1,-0.5) {$n$};

        \graph[nodes={circle, fill}]{
        (a_1) -> (a_2) -> (a_3) -> (a_4);
        };
    }
\]
For $1 \leq i < j \leq n-1$, we denote by $S^i_j$ the unique indecomposable module over $k \vec{\mathbb A}_{n}$ with top $S_i$ and socle $S_j$. We have the following observation.

\begin{lem}
   Let $\{S^{i_1}_{j_1}, \cdots, S^{i_r}_{j_r} \} \in \mathrm{Orth}(k \vec{\mathbb A}_{n} \text{-} \mathrm{mod})$.

   (1) $i_1, \cdots, i_r$ are pairwise distinct.

   (2) Assume that $i_1 < \cdots <i_r$, then $S^{i_1}_{j_1}, \cdots, S^{i_r}_{j_r}$ is an exceptional sequence.
\end{lem}
\begin{proof}
   (1) It is straightforward by the orthogonality and Lemma \ref{uni_hom}.

   (2) We only need to check that $\mathrm{Ext}^1(S^{i_2}_{j_2},S^{i_1}_{j_i})=0$. In fact by Lemma \ref{uni_hom}, $\mathrm{Ext}^1(S^{i_2}_{j_2},S^{i_1}_{j_i})=0$ implies that $i_2>j_1$ or $j_2<j_1$. We now elaborate it.
   
   If $i_2>j_1$, take the projective resolution of $S^{i_2}_{j_2}$
   \[0 \to S^{i_2+1}_{n} \to S^{i_2}_n \to S^{i_2}_{j_2} \to 0,\]
   and apply $\mathrm{Hom}(-,S^{i_1}_{j_1})$ to it, we get an exact sequence
   \[0 \to \mathrm{Hom}(S^{i_2}_{j_2},S^{i_1}_{j_1}) \to \mathrm{Hom}(S^{i_2}_{n},S^{i_1}_{j_1}) \to \mathrm{Hom}(S^{i_2+1}_{n},S^{i_1}_{j_1}) \to \mathrm{Ext}^1(S^{i_2}_{j_2},S^{i_1}_{j_1}) \to 0.\]
   By Lemma \ref{uni_hom}, $\mathrm{Hom}(S^{i_2+1}_{n},S^{i_1}_{j_1})=0$, so $\mathrm{Ext}^1(S^{i_2}_{j_2},S^{i_1}_{j_1}) =0$.

   If $j_2<j_1$, take the injective resolution of $S^{i_1}_{j_1}$
   \[0 \to S^{i_1}_{j_1} \to S^{1}_{j_1} \to S^{1}_{i_1-1} \to 0,\]
   and apply $\mathrm{Hom}(S^{i_2}_{j_2},-)$ to it, we get an exact sequence
   \[0 \to \mathrm{Hom}(S^{i_2}_{j_2},S^{i_1}_{j_1}) \to \mathrm{Hom}(S^{i_2}_{j_2},S^{1}_{j_1}) \to \mathrm{Hom}(S^{i_2}_{j_2},S^{1}_{j_1}) \to \mathrm{Ext}^1(S^{i_2}_{j_2},S^{i_1}_{j_1}) \to 0. \]
   By Lemma \ref{uni_hom}, $\mathrm{Hom}(S^{i_2}_{j_2},S^{1}_{j_1})=0$, so $\mathrm{Ext}^1(S^{i_2}_{j_2},S^{i_1}_{j_1}) =0$.
\end{proof}

Recall that the Weyl group $W_{n-1}$ of $k \vec{\mathbb A}_{n-1}$ is isomorphic to the symmetry group $S_n$ 
\begin{align*}
   W_{n-1} &\xrightarrow{\sim} S_n \\
   s_i &\mapsto (i \ i+1), i=1,\cdots, n-1.
\end{align*}
Through this isomorphism, the reflection $r_{[S^i_j]}$ is send to the transposition $(i,j+1)$: $[S^i_j]=s_j [S^i_{j-1}]$ implies that $r_{[S^i_j]}=s_j r_{[S^i_{j-1}]} s_j^{-1}$ while $j>i$. So Theorem \ref{hubery_krause} for type $A$ can be written as follows. 

\begin{prop} \label{correspondence for type A}
   There is an isomorphism of lattices
   \begin{align*}
      \mathrm{Orth}(k \vec{\mathbb A}_{n-1}) &\xrightarrow{\sim} \mathrm{NC}(S_n, (12\cdots n))\\
      \{S^{i_1}_{j_1}, S^{i_2}_{j_2},\cdots, S^{i_r}_{j_r} \} &\mapsto (i_1 \ j_1+1)(i_2 \ j_2+1) \cdots (i_r \ j_r+1),
   \end{align*}
   where $i_1 <i_2< \cdots <i_r$.
\end{prop}

The proposition above is essentially the same with the statement in \cite[Lemma 6.2]{krause2021category}. They construct a map as follows: let $(a_1< a_2< \cdots< a_p)$ be a block of $\mathrm{NC}^A(n)$, it corresponds to a set of indecomposable objects $S^{a_1}_{a_2-1},S^{a_2}_{a_3-1},\cdots,S^{a_{p-1}}_{a_p-1}$ in $k \vec{\mathbb A}_{n-1} \text{-mod}$. Note that $r_{[S^{a_i}_{a_{i+1}-1}]}$ corresponds to $(a_i a_{i+1})$ in $S_n$, and $(a_1 a_2)(a_2 a_3)\cdots (a_{p-1} a_p)=(a_1 a_2 \cdots a_p)$.

Now we consider bricks and wide subcategories of $\mathcal T_n$. According to \cite[Subsection~1.8]{chen2009introduction}, we can assume that $\{T_1,\cdots,T_n \}$ is a complete set of representatives of the isomorphic classes of simple objects in $\mathcal T_n$, and $\tau T_i=T_{i+1}$, for any $i=1,2, \cdots n$. Here, we set $T_{n+1}:=T_1$. The following description of bricks and exceptional objects in a tube is well-known for experts.

\begin{prop} [{\cite[Corollary 2.2.5]{dichev2009thick}, \cite[Lemma 4.2]{cheng2024wide}}]\label{tube indec length}
   Let $E \in \mathcal T_n$ be an indecomposable object. Then

   (1) $E$ is a brick if and only if $l(E) \leq n$;

   (2) $E$ is exceptional if and only if $l(E)<n$.
\end{prop}

Furthermore, we have the following criterion for whether a wide subcategory is finitely generated.

\begin{prop} \cite[Proposition 4.4(3)]{cheng2024wide} \label{prop nexc wide subcat tube}
   Let $\mathcal W \in \mathrm{Wid}(\mathcal T_n)$. Then $\mathcal W \in \mathrm{NExc}(\mathcal T_n)$ if and only if there exists an indecomposable object $U \in \mathcal W$ such that $l(U)=n$.
\end{prop}

We shows the perpendicular category of a brick in $\mathcal T_n$.

\begin{prop} \label{prop tube brick perpendicular}
   Let $1< m < n-1$, and let $S \in \mathcal T_n$ be a simple object. Then the following hold.

   (1) $\prescript{\perp}{}{S}=\langle \tau S, \tau^2 S, \cdots, \tau^{n-2} S, \tau^{n-1}S^{[2]} \rangle$. In particular, 
   \[
      \prescript{\perp}{}{T_1}=\langle T_2, T_3, \cdots, T_{n-1}, T_n^{[2]} \rangle.
   \]

   (2) $\prescript{\perp}{}{S^{[m]}}=\langle S, \cdots, \tau^{m-2}S,\tau^m S, \cdots, \tau ^{n-2} S, \tau ^{n-1} S^{[m+1]} \rangle$. In particular,
   \[
      \prescript{\perp}{}{T_1^{[m]}}=\langle T_1, \cdots, T_{m-1},T_{m+1}, \cdots, T_{n-1}, T_n^{[m+1]} \rangle.
   \]

   (3) $\prescript{\perp}{}{S^{[n-1]}}=\langle S, \cdots, \tau^{n-3}S, \tau ^{n-1} S^{[n]} \rangle$. In particular,
   \[
      \prescript{\perp}{}{T_1^{[n-1]}}=\langle T_1, \cdots, T_{n-2}, T_n^{[n]} \rangle.
   \]

   (4) $\prescript{\perp}{}{S^{[n]}}=\langle S, \tau S, \cdots , \tau^{n-2}S \rangle$. In particular, 
   \[
      \prescript{\perp}{}{T_1^{[n]}}=\langle T_1, T_2, \cdots, T_{n-1} \rangle.
   \]

   (1') $S^{\perp}=\langle S^{[2]}, \tau ^2 S, \cdots, \tau ^{n-1} S \rangle$. In particular, 
   \[
      T_1^{\perp}=\langle T_1^{[2]}, T_3, \cdots, T_{n} \rangle.
   \]

   (2') $(S^{[m]})^{\perp}=\langle S^{[m+1]}, \tau S, \cdots, \tau ^{m-1} S, \tau ^{m+1} S, \cdots, \tau ^{n-1} S\rangle$. In particular,
   \[
      (T_1^{[m]})^{\perp}=\langle T_1^{[m+1]},T_2 \cdots, T_{m},T_{m+2}, \cdots, T_{n} \rangle.
   \]

   (3') $(S^{[n-1]})^{\perp}=\langle S^{[n]}, \tau S, \cdots, \tau^{n-2} S \rangle$, in particular,
   \[
      (T_1^{[n-1]})^{\perp}=\langle T_1^{[n]}, T_2, \cdots, T_{n-1} \rangle.
   \]

   (4') $(S^{[n]})^{\perp}=\langle \tau S, \tau^2 S, \cdots , \tau^{n-1}S \rangle$. In particular, 
   \[
      (T_1^{[n]})^{\perp}=\langle T_2, T_3, \cdots, T_{n} \rangle.
   \]
\end{prop}

\begin{proof}
   For (1), (2)and (3), we only need to prove that for $1<m<n-1$,
   \[
      \prescript{\perp}{}{T_1^{[m]}}=\langle T_1, \cdots, T_{m-1},T_{m+1}, \cdots, T_{n-1}, T_n^{[m+1]} \rangle.
   \]
   It is straightforward to check that 
   \[
      T_1, \cdots, T_{m-1},T_{m+1}, \cdots, T_{n-1}, T_n^{[m+1]} \in \prescript{\perp}{}{T_1^{[m]}},
   \]
   and $T_1, \cdots, T_{m-1},T_{m+1}, \cdots, T_{n-1}, T_n^{[m+1]}$ is a set of orthogonal bricks. For the other direction $\prescript{\perp}{}{T_1^{[m]}} \subseteq \langle T_1, \cdots, T_{m-1},T_{m+1}, \cdots, T_{n-1}, T_n^{[m+1]} \rangle$, 
   take an indecomposable object $U$ in $\prescript{\perp}{}{T_1^{[m]}}$. Then $\mathrm{Hom}_{\mathcal T_n}(U,T_1^{[m]})=0$ implies that if $T_m$ is a composition factor of $U$, then $T_n^{[m+1]}$ is a factor of $U$. And by Serre duality, $\mathrm{Ext}_{\mathcal T_n}^1 (U,T_1^{[m]})=0$ if and only if $\mathrm{Hom}_{\mathcal T_n} (T_n^{[m]},U)=0$, which implies that if $n$ is a composition factor of $U$, then $T_n^{[m+1]}$ is a factor of $U$. That is 
   \begin{align*}
      \prescript{\perp}{}{T_1^{[m]}} &\subseteq \mathcal U(T_1, \cdots, T_{m-1},T_{m+1}, \cdots, T_{n-1}, T_n^{[m+1]}) \\
      &=\langle T_1, \cdots, T_{m-1},T_{m+1}, \cdots, T_{n-1}, T_n^{[m+1]}\rangle.
   \end{align*}

   For the proof of (4), we refer to \cite[Lemma 4.5(1)]{cheng2024wide}.
\end{proof}

\begin{figure}
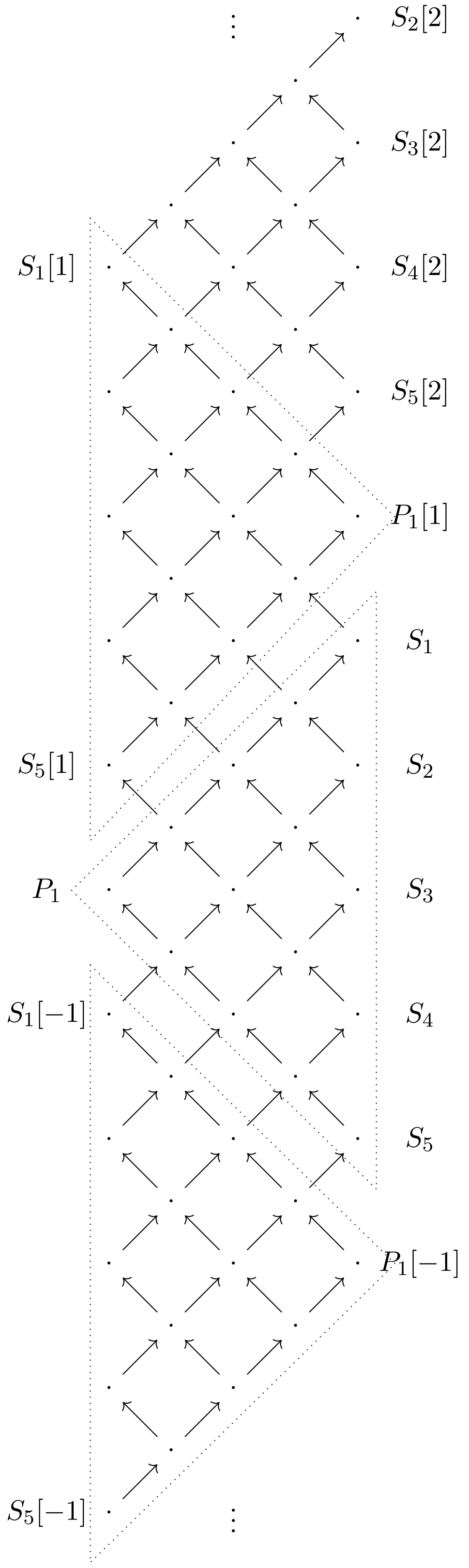

   \tikz[rotate=90,scale=1.7] {
       \node (a1) at (-2,2) {$\cdot$};
       \node (a2) at (-1,2) {$\cdot$};
       \node (a3) at (0,2) {$\cdot$};
       \node (a4) at (1,2) {$\cdot$};
       \node (a5) at (2,2) {$\cdot$};
       \node (a6) at (3,2) {$\cdot$};
       \node (a7) at (4,2) {$\cdot$};
       \node (a8) at (5,2) {$\cdot$};
       \node (a9) at (6,2) {$\cdot$};
       \node (a10) at (7,2) {$\cdot$};
       \node (a11) at (8,2) {$\cdot$};

       \node (b1) at (-1.5,1.5) {$\cdot$};
       \node (b2) at (-0.5,1.5) {$\cdot$};
       \node (b3) at (0.5,1.5) {$\cdot$};
       \node (b4) at (1.5,1.5) {$\cdot$};
       \node (b5) at (2.5,1.5) {$\cdot$};
       \node (b6) at (3.5,1.5) {$\cdot$};
       \node (b7) at (4.5,1.5) {$\cdot$};
       \node (b8) at (5.5,1.5) {$\cdot$};
       \node (b9) at (6.5,1.5) {$\cdot$};
       \node (b10) at (7.5,1.5) {$\cdot$};
       \node (b11) at (8.5,1.5) {$\cdot$};

       \node (c1) at (-1,1) {$\cdot$};
       \node (c2) at (0,1) {$\cdot$};
       \node (c3) at (1,1) {$\cdot$};
       \node (c4) at (2,1) {$\cdot$};
       \node (c5) at (3,1) {$\cdot$};
       \node (c6) at (4,1) {$\cdot$};
       \node (c7) at (5,1) {$\cdot$};
       \node (c8) at (6,1) {$\cdot$};
       \node (c9) at (7,1) {$\cdot$};
       \node (c10) at (8,1) {$\cdot$};
       \node (c11) at (9,1) {$\cdot$};

       \node (d1) at (-0.5,0.5) {$\cdot$};
       \node (d2) at (0.5,0.5) {$\cdot$};
       \node (d3) at (1.5,0.5) {$\cdot$};
       \node (d4) at (2.5,0.5) {$\cdot$};
       \node (d5) at (3.5,0.5) {$\cdot$};
       \node (d6) at (4.5,0.5) {$\cdot$};
       \node (d7) at (5.5,0.5) {$\cdot$};
       \node (d8) at (6.5,0.5) {$\cdot$};
       \node (d9) at (7.5,0.5) {$\cdot$};
       \node (d10) at (8.5,0.5) {$\cdot$};
       \node (d11) at (9.5,0.5) {$\cdot$};

       \node (e1) at (0,0) {$\cdot$};
       \node (e2) at (1,0) {$\cdot$};
       \node (e3) at (2,0) {$\cdot$};
       \node (e4) at (3,0) {$\cdot$};
       \node (e5) at (4,0) {$\cdot$};
       \node (e6) at (5,0) {$\cdot$};
       \node (e7) at (6,0) {$\cdot$};
       \node (e8) at (7,0) {$\cdot$};
       \node (e9) at (8,0) {$\cdot$};
       \node (e10) at (9,0) {$\cdot$};
       \node (e11) at (10,0) {$\cdot$};

       \node(l) at (-2,1) {$\vdots$};
       \node(r) at (10,1) {$\vdots$};

       \node (f1) at (0,-0.5) {$P_1[-1]$};
       \node (f2) at (1,-0.5) {$S_5$};
       \node (f3) at (2,-0.5) {$S_4$};
       \node (f4) at (3,-0.5) {$S_3$};
       \node (f5) at (4,-0.5) {$S_2$};
       \node (f6) at (5,-0.5) {$S_1$};
       \node (f7) at (6,-0.5) {$P_1[1]$};
       \node (f8) at (7,-0.5) {$S_5[2]$};
       \node (f8) at (8,-0.5) {$S_4[2]$};
       \node (f10) at (9,-0.5) {$S_3[2]$};
       \node (f11) at (10,-0.5) {$S_2[2]$};

       \node (g1) at (-2,2.5) {$S_5[-1]$};
       \node (g5) at (2,2.5) {$S_1[-1]$};
       \node (g6) at (3,2.5) {$P_1$};
       \node (g7) at (4,2.5) {$S_5[1]$};
       \node (g11) at (8,2.5) {$S_1[1]$};
   
       \graph{
           (a1)->(b1)->(c1)->(d1)->(e1);
           (a2)->(b2)->(c2)->(d2)->(e2);
           (a3)->(b3)->(c3)->(d3)->(e3);
           (a4)->(b4)->(c4)->(d4)->(e4);
           (a5)->(b5)->(c5)->(d5)->(e5);
           (a6)->(b6)->(c6)->(d6)->(e6);
           (a7)->(b7)->(c7)->(d7)->(e7);
           (a8)->(b8)->(c8)->(d8)->(e8);
           (a9)->(b9)->(c9)->(d9)->(e9);
           (a10)->(b10)->(c10)->(d10)->(e10);
           (a11)->(b11)->(c11)->(d11)->(e11);
           (b1)->(a2);
           (c1)->(b2)->(a3);
           (d1)->(c2)->(b3)->(a4);
           (e1)->(d2)->(c3)->(b4)->(a5);
           (e2)->(d3)->(c4)->(b5)->(a6);
           (e3)->(d4)->(c5)->(b6)->(a7);
           (e4)->(d5)->(c6)->(b7)->(a8);
           (e5)->(d6)->(c7)->(b8)->(a9);
           (e6)->(d7)->(c8)->(b9)->(a10);
           (e7)->(d8)->(c9)->(b10)->(a11);
           (e8)->(d9)->(c10)->(b11);
           (e9)->(d10)->(c11);
           (e10)->(d11);
       };

       \draw[dotted] (0.6,-0.15) -- (3,2.3) -- (5.4,-0.15) -- (0.6,-0.15);
       \draw[dotted] (0,-0.3) -- (2.4,2.15) -- (-2.4,2.15) -- (0,-0.3);
       \draw[dotted] (6,-0.3) -- (8.4,2.15) -- (3.4,2.15) -- (6,-0.3);
   }
   \caption{The AR quiver of $\mathcal{D}^b(k \vec{\mathbb A}_5 \text{-mod})$}
   \label{AR derived A_5}
\end{figure}

Based on the construction in Proposition~\ref{correspondence for type A}, \cite[Lemma 8.2, Corollary 9.2(2)]{krause2021category} also shows that $\mathrm{Wid}(\mathcal T_n) \cong \mathrm{NC}^B(n)$. We provide an aspect from derived category, and prove that this isomorphism preserves Kreweras complement. First, let us recall the Auslander-Reiten quiver of $\mathcal{D}^b(k \vec{\mathbb A}_n \text{-mod})$. Each indecomposable object in $\mathcal{D}^b(k \vec{\mathbb A}_n \text{-mod})$ is isomorphic to $E[i]$, where $E$ is an indecomposable object in $k \vec{\mathbb A}_n \text{-mod}$ and $i \in \mathbb Z$. Note that $\tau ^{2n} \cong [-2]$, where $\tau$ is the AR translation. Figure \ref{AR derived A_5} shows the AR quiver of $\mathcal{D}^b(k \vec{\mathbb A}_5 \text{-mod})$.

The AR translation $\tau$ acts on $\mathrm{Thick}(\mathcal D^b(k \vec{\mathbb A}_{2n-1} \text{-} \mathrm{mod}))$ naturally. The lattice formed by $\tau ^n$-invariant thick subcategories is isomorphic to $\mathrm{NC}^B(n)$.

\begin{prop} \label{prop tau_n inv Db(A_{2n-1}-mod)}
   There is an isomorphism of lattices
   \[
      \mathrm{Thick}(\mathcal D^b(k \vec{\mathbb A}_{2n-1} \text{-} \mathrm{mod}))^{\tau ^n} \cong \mathrm{NC}^B(n).
   \]
\end{prop}

\begin{proof}
   First, $\mathrm{Thick}(\mathcal D^b(k \vec{\mathbb A}_{2n-1} \text{-} \mathrm{mod})) \cong \mathrm{NC}^A(2n)$. From the AR-quiver of $\mathcal D^b(k \vec{\mathbb A}_{2n-1} \text{-} \mathrm{mod})$, it is clear that the AR translation $\tau$ acts on $\mathrm{NC}^A(2n)$ naturally by $1/2n$ of a full clockwise rotation of the circular representation. So $\mathrm{Thick}(\mathcal D^b(k \vec{\mathbb A}_{2n-1} \text{-} \mathrm{mod}))^{\tau ^n}$ is isomorphic to the lattice $\mathrm{NC}^A(2n)^a$, where $a$ is the antipodal map, i.e. $\mathrm{NC}^B(n)$. 
\end{proof}

We construct an isomorphism from the lattice of orthogonal bricks of $\mathcal T_n$ to the lattice of $\tau ^n$-invariant thick subcategories of $\mathcal D^b(k \vec{\mathbb A}_{2n-1} \text{-} \mathrm{mod})$.

\begin{prop} \label{prop orth Tn to Thick Db}
   There is an isomorphism of lattices
   \begin{align*}
      \mathrm{Orth}(\mathcal T_n) &\xrightarrow{\sim} \mathrm{Thick}(\mathcal D^b(k \vec{\mathbb A}_{2n-1} \text{-} \mathrm{mod}))^{\tau ^n} \\
      \{T_{i_1}^{[j_1]},\cdots,T_{i_r}^{[j_r]} \} &\mapsto \langle S_{i_1}^{[j_1]},\cdots,S_{i_r}^{[j_r]},\tau ^n S_{i_1}^{[j_1]},\cdots, \tau ^n S_{i_r}^{[j_r]}\rangle.
   \end{align*}
\end{prop}

\begin{proof}
   Recall that in Proposition \ref{correspondence for type A}, the inverse of the lattice isomorphism $\mathrm{Orth}(k \vec{\mathbb A} \text{-} \mathrm{mod}) \xrightarrow{\sim} \mathrm{NC}^A(2n)$ is clear: a block $\{a_1<a_2< \cdots <a_p \}$ corresponds to a set of indecomposable objects $S^{a_1}_{a_2-1},S^{a_2}_{a_3-1},\cdots,S^{a_{p-1}}_{a_p-1}$ in $k \vec{\mathbb A}_{2n-1} \text{-mod}$. Recall Proposition~\ref{prop tau_n inv Db(A_{2n-1}-mod)} that there is an isomorphism of lattices 
   \[
      \mathrm{Thick}(\mathcal D^b(k \vec{\mathbb A}_{2n-1} \text{-} \mathrm{mod}))^{\tau ^n} \cong \mathrm{NC}^B(n),
   \]
   so by Corollary \ref{cor brick wide cat}, each $\tau^n$-invariant thick subcategory of $\mathcal D^b(k \vec{\mathbb A}_{2n-1} \text{-} \mathrm{mod})$ is uniquely determined by a set $\mathcal S$ of orthogonal bricks in $k \vec{\mathbb A}_{2n-1} \text{-} \mathrm{mod}$ such that each brick in $\mathcal S$ has top in $\{S_1,S_2,\cdots S_n \}$ and has length not greater than $n$. Precisely, the thick subcategory equals to $\langle \mathcal S, \tau^n \mathcal S \rangle$. Suppose that $\mathcal S=S_{i_1}^{[j_1]},\cdots,S_{i_r}^{[j_r]}$, where $1\leq i_1,\cdots,i_r\leq n$, $j_1,\cdots, j_r \leq n$. Then $\mathcal S$ corresponds to a unique orthogonal set $\mathcal S'=T_{i_1}^{[j_1]},\cdots,T_{i_r}^{[j_r]}$ of bricks in $\mathcal T_n$. In fact, by Proposition \ref{tube indec length}(1), each object in $\mathcal S'$ is a brick in $\mathcal T_n$. And by Lemma \ref{uni_hom}, $T_{i_u}^{[j_u]}$ and $T_{i_v}^{[j_v]}$ are orthogonal if and only if $S_{i_u}^{[j_u]}$ and $S_{i_v}^{[j_v]}$ are orthogonal.
\end{proof}

The following theorem gives a proof of \cite[corollary~9.2~(2)]{krause2021category} using derived category. For another representation-theoretic proof, we refer to \cite[Corollary~11.3]{hanson2025para}. Additionally, we prove that taking left perpendicular category in $\mathrm{Wid}(\mathcal T_n)$ corresponds to the Kreweras complement in $\mathrm{NC}^n(B)$.

\begin{thm} \label{thm tube B}
   There is a natural isomorphism
   \[
   \phi: \mathrm{Wid}(\mathcal T_n) \cong \mathrm{NC}^B(n).
   \]
   Moreover, taking left perpendicular category corresponds to the Kreweras complement via this isomorphism, i.e. for any wide subcategory $\mathcal W$ of $\mathcal T_n$,
   \[
   \phi(\prescript{\perp}{}{\mathcal W})=K(\phi(\mathcal W)).
   \]
\end{thm}

\begin{proof}
   By Corollary~\ref{cor brick wide cat}, Proposition~\ref{prop orth Tn to Thick Db} and Proposition~\ref{prop tau_n inv Db(A_{2n-1}-mod)}, we have
   \[\phi: \mathrm{Wid}(\mathcal T_n) \xrightarrow{\sim} \mathrm{Orth}(\mathcal T_n) \xrightarrow{\sim} \mathrm{Thick}(\mathcal D^b(k \vec{\mathbb A}_{2n-1} \text{-} \mathrm{mod}))^{\tau ^n} \xrightarrow{\sim} \mathrm{NC}^B(n).\]
   By Proposition~\ref{Kreweras in cat}, the isomorphism 
   \[\phi':\mathrm{Thick}(\mathcal D^b(k \vec{\mathbb A}_{2n-1} \text{-} \mathrm{mod}))^{\tau ^n} \xrightarrow{\sim} \mathrm{NC}^B(n)\]
   satisfies $\phi'(\prescript{\perp}{}{\mathcal T})=K(\phi'(\mathcal T))$ for any $\mathcal T \in \mathrm{Thick}(\mathcal D^b(k \vec{\mathbb A}_{2n-1} \text{-} \mathrm{mod}))^{\tau ^n}$. So we only need to prove that the isomorphism
   \[
   \phi'':\mathrm{Wid}(\mathcal T_n) \cong \mathrm{Thick}(\mathcal D^b(k \vec{\mathbb A}_{2n-1} \text{-} \mathrm{mod}))^{\tau ^n}
   \]
   preserves taking left perpendicular categories, i.e. $\phi''(\prescript{\perp}{}{\mathcal W})=\prescript{\perp}{}{\phi''(\mathcal W)}$ for any $\mathcal W \in \mathrm{Wid}(\mathcal T_n)$. 
   
   We prove this by induction on the number $m=|\mathcal S(\mathcal W)|$ for $\mathcal W \in \mathrm{Wid}(\mathcal T_n)$. Recall that $\mathcal S:\mathrm{Wid}(\mathcal T_n) \cong \mathrm{Orth}(\mathcal T_n)$ is a bijection. When $m=1$, we only need to prove that for a brick $T_i^{[j]} \in \mathcal T_n$, we have $\phi''(\prescript{\perp}{}{T_i^{[j]}})=\prescript{\perp}{}{\phi''(T_i^{[j]})}$. For the convenience, we assume that $i=1$. By Proposition~\ref{prop tube brick perpendicular}, we get the orthogonal generating bricks of $\prescript{\perp}{}{\phi''(T_1^{[j]})}$: for $1 \leq j \leq n-2$, 
   \[\prescript{\perp}{}{T_1^{[j]}}=\langle T_1, \cdots, T_{j-1}, T_{j+1},\cdots, T_{n-1}, T_n^{[j+1]} \rangle,\]
   and 
   \[\prescript{\perp}{}{T_1^{[n-1]}}=\langle T_1, \cdots, T_{n-2}, T_n^{[n]} \rangle,\]
   and
   \[\prescript{\perp}{}{T_1^{[n]}}=\langle T_1, \cdots, T_{n-1} \rangle.\]
   So for $1 \leq j \leq n-2$, denote by 
   \[\mathcal U=\{ S_1, \cdots, S_{j-1}, S_{j+1},\cdots, S_{n-1}, S_n^{[j+1]}\} \subseteq \mathcal D^b(k \vec{\mathbb A}_{2n-1} \text{-} \mathrm{mod}), \]
   then
   \begin{align*}
      \phi''(\prescript{\perp}{}{T_1^{[j]}})&=\langle \mathcal U, \tau ^n \mathcal U\rangle, 
   \end{align*}
   with orthogonal generating bricks
   \begin{align*}
      S_1, \cdots, S_{j-1}, &S_{j+1},\cdots, S_{n-1}, S_n^{[j+1]},\\S_{n+1},\cdots, S_{n-j+1}, &S_{n+j+1}, \cdots, S_{2n-1}, S_{j+1}^{[2n-j-1]}
   \end{align*}
   in $k \vec{\mathbb A}_{2n-1} \text{-} \mathrm{mod}$. On the other hand, 
   $\prescript{\perp}{}{\phi''(T_1^{[j]})}$ corresponds to the noncrossing partition
   \begin{align*}
      &K(\{\{1,j+1\},\{n+1,n+j+1\}\})\\
      &=\{ \{1,\cdots,m \},\{m+1,\cdots,n,n+j+1,\cdots, 2n\},\{n+1,\cdots,n+j\} \}
   \end{align*}
   in $\mathrm{NC}^{A}(2n)$, so has the same orthogonal generating bricks in $k \vec{\mathbb A}_{2n-1} \text{-} \mathrm{mod}$ with $\phi''(\prescript{\perp}{}{T_1^{[j]}})$. That is 
   \[
      \phi''(\prescript{\perp}{}{T_i^{[j]}})=\prescript{\perp}{}{\phi''(T_i^{[j]})} \text{ for } 1 \leq j \leq n-2.
   \]
   For $j=n-1$ and $j=n$, the proof is similar and easier. So for $\mathcal W \in \mathrm{Wid}(\mathcal T_n)$ with $\mathcal S(\mathcal W)=1$, we have checked that $\phi''(\prescript{\perp}{}{\mathcal W})=\prescript{\perp}{}{\phi''(\mathcal W)}$.

   Suppose that for $\mathcal W \in \mathrm{Wid}(\mathcal T_n)$ with $\mathcal S(\mathcal W) \leq M$, we have  
   \[
      \phi''(\prescript{\perp}{}{\mathcal W})=\prescript{\perp}{}{\phi''(\mathcal W)},
   \]
   then for $\mathcal W$ with $\mathcal S(\mathcal W) = M$, there exists $\mathcal U, \mathcal V \in \mathrm{Wid}(\mathcal T_n)$ such that 
   \[
      \mathcal W=\mathcal U+\mathcal V \text{ and } \mathcal S(\mathcal U)<\mathcal S(\mathcal W), \ \mathcal S(\mathcal V)<\mathcal S(\mathcal W).
   \]
   Note that $\phi''$ is an isomorphism of lattices, so it preserves joins and meets. From Proposition~\ref{tube perpendicular} we know that taking perpendicular category gives duality for the lattice $\mathrm{Wid}(\mathcal T_n)$, and the same to the lattice $\mathrm{Wid}(k \vec{\mathbb A} \text{-} \mathrm{mod})$ from Theorem~\ref{hubery_krause}. So
   \begin{align*}
      \prescript{\perp}{}{\phi''(\mathcal W)}&=\prescript{\perp}{}{\phi''(\mathcal U+\mathcal V)}=\prescript{\perp}{}{\phi''(\mathcal U)} \cap \prescript{\perp}{}{\phi''(\mathcal V)} \\
      &=\phi''(\prescript{\perp}{}{\mathcal U}) \cap \phi''(\prescript{\perp}{}{\mathcal V}) =\phi''(\prescript{\perp}{}{\mathcal U} \cap \prescript{\perp}{}{\mathcal V}) \\
      &=\phi(\prescript{\perp}{}{(\mathcal U+\mathcal V)})=\phi''(\prescript{\perp}{}{\mathcal W}).
   \end{align*}
   The proof is complete.
\end{proof}

\begin{prop} \label{prop subcat zero block}
   Let $\mathcal W \in \mathrm{Wid}(\mathcal T_n)$, then $\mathcal W \in \mathrm{NExc}(\mathcal T_n)$ if and only if $\phi(\mathcal W) \in \mathrm{NC}^B(n)$ has a zero block, where $\phi:\mathrm{Wid}(\mathcal T_n) \cong \mathrm{NC}^B(n)$ is defined in Theorem~\ref{thm tube B}.
\end{prop}

\begin{proof}
   $\Rightarrow$: By Proposition~\ref{prop nexc wide subcat tube}, $\mathcal W \in \mathrm{NExc}(\mathcal T_n)$ if and only if $\mathcal W$ contains an indecomposable object $U$ which has length $n$. Then the composition factors of $U$ in in $\mathcal W$ corresponds to a zero block in $\mathrm{NC}^B(n)$.

   $\Leftarrow$: If $\phi (\mathcal W)$ has a zero block 
   \[\{i_1, i_2,\cdots,i_r, \bar{i}_1, \bar{i}_2, \cdots, \bar{i}_r \}, 1 \leq i_1<i_2<\cdots<i_r \leq n,\]
   then by the definition of $\phi$, we have that 
   \[T_{i_1}^{[i_2-i_1]}, T_{i_2}^{[i_3-i_2]}, \cdots, T_{i_{r-1}}^{[i_r-i_{r-1}]}, T_{i_r}^{[n+i_1-i_r]} \] 
   lies in $\mathcal W$ as simple objects. So $T_{i_1}^{[n]} \in \mathcal W$ as their extensions, and has length $n$ in $\mathcal T_n$. By Proposition~\ref{prop nexc wide subcat tube}, $\mathcal W \in \mathrm{NExc}(\mathcal T_n)$.
\end{proof}

We reobtain a known result by translating Proposition~\ref{tube perpendicular} to a combinatorial version. For a pure combinatorial proof, we refer to \cite[section~1.2]{biane2003non}.

\begin{cor}
   The Kreweras complement in $\mathrm{NC}^B(n)$ gives a duality between noncrossing partitions which has a zero block and those has not zero blocks.
\end{cor}

\begin{proof}
   By Theorem~\ref{thm tube B}, the Kreweras complement in $\mathrm{NC}^B(n)$ corresponds to taking left perpendicular categories in $\mathrm{Wid}(\mathcal T_n)$. Then it is straightforward by Proposition~\ref{tube perpendicular} and Proposition~\ref{prop subcat zero block}.
\end{proof}

\section{Weighted projective line}

Let $\bm p=(p_1,p_2,\cdots,p_n)\in \mathbb N_+^n$. There is an abelian group
\[
L(\bm p):=\langle\vec{x_1},\vec{x_2},\cdots,\vec{x_n} \ | \ p_1\vec{x_1}=p_2\vec{x_2}=\cdots=p_n\vec{x_n}\rangle.
\]
The \emph{canonical element} $\vec{c}=p_1\vec{x_1}=\cdots =p_n \vec{x_n}$, and the \emph{dualizing element} $\vec{\omega}=(n-2)\vec{c}-\sum_{i=1}^n\vec{x_i}$. We have the \emph{degree map}
\[
\delta: L(\bm p)\twoheadrightarrow \mathbb Z, \ \vec{x_i}\mapsto p/p_i,
\]
where $p=\mathrm{l.c.m.}(p_1,\cdots,p_n)$. Note that each $\vec l \in L(\bm p)$ can be uniquely written in \emph{normal form}
\[
   \vec l =\sum_{k=1}^n l_k \vec{x_k}+l \vec c \ \text{  with  }\  0 \leq l_i <p_i \ \text{  and  }\  l \in \mathbb Z.
\]

Let $k$ be an algebraically closed field, and let $\bm \lambda=(\lambda_1, \lambda_2, \cdots, \lambda_n)$ be an $n$-tuple of different (closed) points in $\mathbb P^{1}_k$, normalized such that $\lambda_1=\infty, \lambda_2=0, \lambda_3=1$. Let
\[
S=S(\bm p, \bm{\lambda}):=k[X_1,X_2,\cdots,X_n]/(X_i^{p_i}-X_2^{p_2}+\lambda_i X_1^{p_1})_{3 \leq i \leq n}
\]
be an $L(\bm p)$-graded ring with $\mathrm{deg}(X_i)=\vec{x_i}$. We denote by $x_i$ the quotient class of $X_i$ in $S(\bm p, \bm \lambda)$. It is an $L(\bm p)$-graded unique factorization domain; see \cite[Proposition~1.3]{geigle1987class}. Let
\[
   \mathbb X=\mathbb X(\bm p, \bm \lambda):=\{\mathfrak{p} \text{ homogeneous prime ideal } | \ \mathfrak{p} \nsupseteq (x_1, \cdots x_n)\}
\]
be a topological space equipped with the \emph{Zariski topology}. A topological basis of $\mathbb X$ is
\[
   \{D(f):=\{\mathfrak{p} \in \mathbb X \ | \ f \notin \mathfrak{p}\} \ | \ f \in (x_1, \cdots , x_n) \text{ homogeneous}\}.
\]
One defines the \emph{structure sheaf} $\mathcal O_{\mathbb X}$ on $\mathbb X$ by setting  
\[
   \mathcal O_{\mathbb X}(D(f))=S_f:=\{ \frac{s}{f^n} \ | \ s \in S, n \in \mathbb N \}.
\]
Then $\mathcal O_{\mathbb X}$ is a sheaf of $L(\bm p)$-graded rings. The \emph{weighted projective line} is by definition the $L(\bm p)$-\emph{graded locally ringed space} $(\mathbb X, \mathcal O_{\mathbb X})$. 

By \cite[Proposition~1.3]{geigle1987class}, as a set,
\[
   \mathbb X=\{(0),(x_i),(x_2^{p_2}- \lambda x_1^{p_1}) \ | \ i=1,2, \cdots, n, \lambda \notin \bm \lambda \}.
\]
We may view the closed points $(x_1),(x_2), \cdots, (x_n)$ in $\mathbb X$ as points $ \lambda_1, \lambda_2, \cdots, \lambda_n $ in $\mathbb P^1(k)$, respectively. They are called \emph{exceptional points}. And we also view the point $(x_2^{p_2}- \lambda x_1^{p_1}) \in \mathbb X$ as the point $\lambda \in \mathbb P^1(k)$, where $\lambda \neq \lambda_1, \cdots, \lambda_n$. They are called \emph{ordinary points}. Note that the only non-closed point in $\mathbb X$ (also for $\mathbb P^1(k)$) is the generic point $(0)$.

The group $L(\bm{p})$ acts on the the category of $L(\bm p)$-graded $S$-modules by \emph{grading shift}: for an $L(\bm p)$-graded module $M$, and $\vec l \in L(\bm p)$ we denote by $M(\vec l)$ the $S$-module $M$ with the new grading $M(\vec l)_{\vec l'}=M_{\vec l + \vec l'}$. It induces the grading shift functor on the category of \emph{$L(\bm p)$-graded $\mathcal O_{\mathbb X}$-module}. Precisely, for an $\mathcal O_{\mathbb X}$-module $F$, and $\vec l \in L(p)$, one defines $F(\vec l)$ by setting 
\[
   (F(\vec l))(U)=(F(U))(\vec l)
\]
for any open subset $U$ of $\mathbb X$.

A \emph{coherent sheaf} on $\mathbb X$ is by definition an $L(\bm p)$-graded $\mathcal O_{\mathbb X}$-module $F$, where for each $x \in \mathbb X$ there is an open neighborhood $U$ of $x$ and an exact sequence
\[
   \bigoplus_{j=1}^{m'} \mathcal O_{\mathbb X}(\vec l_j)|_U \to \bigoplus_{i=1}^m \mathcal O_{\mathbb X}(\vec l_i)|_U \to F|_U \to 0.
\]

By \cite{geigle1987class}, the category $\mathrm{coh} \text{-} \mathbb X$ is an Ext-finite hereditary abelian $k$-category with Serre duality
\[
\mathrm{Ext}_{\mathbb X}^1(F,G)\xrightarrow{\sim} D \mathrm{Hom}_{\mathbb X}(G,F(\vec \omega)),
\]
where $D$ denotes the dual functor $\mathrm{Hom}_k (-,k)$ of $k$-vector spaces.

A coherent sheaf $F$ is called \emph{locally free} (or a \emph{vector bundle}), if for any $x\in \mathbb X$, there exists an open neighborhood $U\ni x$, $r\in \mathbb N$, and $\vec l_1,\cdots ,\vec l_r \in L(\bm p)$, such that $F|_{U}\cong \bigoplus_{i=1}^r\mathcal O_{\mathbb X}(\vec l_i)|_{U}$. We call $r$ the \emph{rank} of $F$, denoted by $\mathrm{rk}(F)$. We write $\mathrm{vect}(\mathbb X)$ for the subcategory of vector bundles. Rank one vector bundles are called \emph{line bundles}. Each line bundle has the form $\mathcal O_{\mathbb{X}}(\vec l)$, where $\vec l \in L(\bm p)$. Each rank $n$ vector bundle $F$ has a filtration
\[
0=F_0 \subset F_1 \subset \cdots \subset F_n=F,
\]
whose factors $F_i/F_{i-1}$ are line bundles.

Each $F \in \mathrm{coh} \text{-} \mathbb X$ has a unique maximal finite length subsheaf, denoted by $\mathrm{t}F$. For any $F \in \mathrm{coh} \text{-} \mathbb X$, the short exact sequence
\[
0 \to \mathrm{t}F \to F \to F/\mathrm{t}F \to 0
\]
splits, and $F/\mathrm{t}F$ is a vector bundle.

We denote by $\mathrm{coh}_0 \text{-} \mathbb X$ the subcategory consisting of all \emph{finite length sheaves} (i.e. \emph{torsion sheaves}). A coherent sheaf $F$ is torsion if and only if its \emph{support} 
\[
   \mathrm{supp}(F):=\{x \in \mathbb X \ | \ F_x \neq 0 \}
\]
is a finite set. And in this case, we have $F=\oplus_{x \in \mathrm{supp}(F)} F_x$.

We denote by
\[
   \mathrm{tor}_{\lambda} \text{-} \mathbb X:=\{ F \in \mathrm{coh}_0 \text{-} \mathbb X \ | \ \mathrm{supp}(F)=\{ \lambda \} \}
\]
the category of torsion sheaves concentrated on the point $\lambda$. According to \cite[Subsection~2.4]{geigle1987class}, for $\lambda \in \mathbb X$ of weight $w$, $\mathrm{tor}_{\lambda} \text{-} \mathbb X$ is equivalent to the $n$-tube $\mathcal T_n$. It is a wide subcategory of $\mathrm{coh} \text{-} \mathbb X$.

Let $u_\lambda=x_2^{p_2}-\lambda x_1^{p_1} \in S(\bm p, \bm \lambda)$ with $\lambda \neq \lambda_1, \cdots, \lambda_n$. Then each simple sheaf on $\mathbb X$ is isomorphic to one of those $S_\lambda, S_{i,j}$ defined as the third term of the following exact sequences:
\[
 0 \to \mathcal O_{\mathbb X}\xrightarrow{u_\lambda} \mathcal O_{\mathbb X}(\vec c)\to S_\lambda \to 0,
\]
\[
 0 \to \mathcal O_{\mathbb X}((j-1)\vec {x_i})\xrightarrow{x_i} \mathcal O_{\mathbb X}(j\vec{x_i})\to S_{i,j} \to 0,
\]
where $u_\lambda$ and $x_i$ denote the morphisms induced by multiplications with $u_\lambda$ and $x_i$, respectively; $i=1,2, \cdots, n$; $j=1,2,\cdots,p_i \in \mathbb Z/ p_i \mathbb Z$. Note that for any $\vec l =\sum_{k=1}^n l_k \vec{x_k}+l \vec c \in L(\bm p)$ written in normal form, we have 
\[
   S_\lambda (\vec l) \cong S_\lambda, S_{i,j}(\vec l) \cong S_{i,j+l_k}.
\]
In particular, we have 
\[
   S_{\lambda}(\vec \omega) \cong S_\lambda, S_{i,j}(\vec \omega) \cong S_{i, j-1}.
\]

Following \cite{cheng2024wide}, we use the notion of $\vec c$-invariant wide subcategory to study wide subcategories of a weighted projective line. A wide subcategory $\mathcal W$ is called \emph{$\vec c$-invariant} if $\mathcal W(\vec c)=\mathcal W$, where 
\[\mathcal W(\vec c):=\{F \in \mathrm{coh}\text{-}\mathbb X \ | \ F(\vec c)=F\}.\]
Note that if $\mathcal W$ is $\vec c$-invariant, then $\mathcal W(-\vec c)=\mathcal W$ since $\mathcal W$ is wide. And for a torsion sheaf $T$, we have $T \cong T(\vec c)$. So a wide subcategory formed by torsion sheaves is always $\vec c$-invariant. We denote by $\mathrm{Wid}(\mathrm{coh} \text{-} \mathbb X)^{\vec c}$ the set of $\vec c$-invariant wide subcategories, and $\mathrm{Exc}(\mathrm{coh} \text{-} \mathbb X)^{\vec c}$ the set of finitely generated $\vec c$-invariant subcategories.

The next proposition states that a non-exceptional indecomposable torsion sheaf will force a wide subcategory containing it to be $\vec c$-invariant.

\begin{prop} [{\cite[Proposition 6.2]{cheng2024wide}}] \label{contain_non_exc}
Let $\mathcal W$ be a wide subcategory of $\mathrm{coh} \text{-} \mathbb X$ containing a vector bundle. Then $\mathcal W$ contains a non-exceptional indecomposable torsion sheaf if and only if $\mathcal W$ is $\vec c$-invariant.
\end{prop}

The following is the structure theorem of $\vec c$-invariant wide subcategories.

\begin{thm} [{\cite[Theorem 6.3]{cheng2024wide}}]\label{c_inv}
   Let $\mathbb X$ be a weighted projective line with $n$ exceptional points, and let $\mathcal W$ be a $\vec c$ -invariant wide subcategory of $\mathrm{coh} \text{-} \mathbb X$ containing a vector bundle. Then
   
   (1) For any $i=1,2,\cdots, n$, there exists an (unnecessarily non-zero) exceptional sequence $E^i=(E_1^i,\cdots,E_{r_i}^i)$ in $\mathrm{tor}_{\lambda_i} \text{-} \mathbb X$, such that
   \[
      \mathcal W \cap \mathrm{tor}_{\lambda_i} \text{-} \mathbb X =(E^i)^{\perp} \cap \mathrm{tor}_{\lambda_i} \text{-} \mathbb X.
   \]
   
   (2) $\mathcal W=(E^1,\cdots,E^n)^{\perp}$, where $E^1, \cdots , E^n$ are defined in (1).
\end{thm}

The following two well-known propositions state the structures of perpendicular categories in $\mathrm{coh} \text{-} \mathbb X$. We provide a proof for Proposition~\ref{prop bundle perpendicular} for the difficulty to obtain literatures \cite{huebner1993categories} and \cite{huebner1996thesis}.

\begin{prop} [{\cite{huebner1993categories}, \cite[Theorem~5.4]{huebner1996thesis}}] \label{prop bundle perpendicular}
   Let $\mathbb X$ be a weighted projective line. $E \in \mathrm{coh} \text{-} \mathbb X$ be a rigid vector bundle. Then $E^{\perp}$ is equivalent to $A \text{-} \mathrm{mod}$ for some finite dimensional hereditary algebra $A$.
\end{prop}

\begin{proof}
   Since each vector bundle $F$ has a filtration of line bundles, for any line bundle $L$, there exists $\vec l \in L(\bm p)$, such that for any $\vec l' \geq \vec l$, $L(- \vec l')$ can be imbedded to $F$. So for the canonical tilting object $T=\oplus_{0 \leq l \leq \vec c} \mathcal O(\vec l)$, there is a sufficiently large $\vec l \in L(\bm p)$ such that $T(- \vec l)$ can be cogenerated by $E$. Notice that $T(- \vec l)$ is also a tilting object, by Proposition~\ref{prop perp cat is module cat} , this proposition holds.
\end{proof}

\begin{prop} [{\cite[Proposition~6.5]{cheng2024wide}}] \label{prop tor perpendicular}
   Let $\mathbb X(p_1, \cdots, p_n)$ be a weighted projective line, $E \in \mathrm{coh} \text{-} \mathbb X$ a torsion exceptional sheaf concentrated on the point $\lambda_i$. Then 
   \[E^{\perp} \cong \mathrm{coh} \text{-} \mathbb X' \coprod k \vec{\mathbb A}_{l(E)-1} \text{-} \mathrm{mod},     \]
   where $\mathbb X'$ is of weight type $(p_1, \cdots, p_{i}-l(E), \cdots, p_n)$.
\end{prop}

Weighted projective lines are classified by the value of $\delta(\vec \omega)$:

(1) If $\delta(\vec \omega)<0$, then $\mathbb X$ is called of \emph{domestic} type. Precisely, $\bm{p}=(p,q)$, $(2,2,n)$, $(2,3,3)$, $(2,3,4)$, $(2,3,5)$.

(2) If $\delta(\vec \omega)=0$, then $\mathbb X$ is called of \emph{tubular} type. Precisely, $\bm{p}=(2,2,2,2)$, $(3,3,3)$, $(2,4,4)$, $(2,3,6)$.

(3) If $\delta(\vec \omega)>0$, then $\mathbb X$ is called of \emph{wild} type.

The classification of indecomposable vector bundles and Auslander-Reiten quiver for a domestic weighted projective line is well-known, see the following theorem.

\begin{thm}  [{\cite[Theorem~3.5]{lenzing2011weighted}}]\label{domestic}
   Let $\mathbb X$ be a weighted projective line of domestic type, then the following properties hold: 
   
   (1) Each indecomposable bundle is exceptional. 
   
   (2) Let $T$ be the direct sum of all indecomposable bundles $E$ satisfying $0 \leq \mu(E) < -\delta(\omega)$. Then $T$ is a tilting object.
   
   (3) The endomorphism ring of $T$ is isomorphic to the path algebra of an extended Dynkin quiver, which is either of type $\tilde{\mathbb A}_{p,q}$ or of type $\tilde{\mathbb D}_n$ ($n \geq 4$), $\tilde{\mathbb E}_6$, $\tilde{\mathbb E}_7$ or $\tilde{\mathbb E}_8$ with bipartite orientation.
\end{thm}

We have the following classification theorem for wide subcategories over a weighted projective line of domestic type.

\begin{prop} [{\cite[Proposition 7.1]{cheng2024wide}}] \label{thick_domestic_bundle}
   Let $\mathbb X$ be a weighted projective line of domestic type and $\mathcal W$ a wide subcategory of $\mathrm{coh} \text{-} \mathbb X$. If $\mathcal W$ contains a vector bundle, then $\mathcal W \in \mathrm{Exc}(\mathrm{coh} \text{-} \mathbb X)$.
\end{prop}

The following proposition is known for experts.

\begin{prop} \label{prop domestic rep-finite}
   Let $\mathbb X$ be a weighted projective line of domestic type, $E$ be an exceptional vector bundle. Then there exists a finite dimensional representation-finite hereditary algebra $A$ such that $E^{\perp} \cong A \text{-} \mathrm{mod}$.
\end{prop}

\begin{proof}
   By Proposition~\ref{prop bundle perpendicular}, There is a finite dimensional hereditary algebra $A$ such that $E^{\perp} \cong A \text{-} \mathrm{mod}$. We only need to prove that $A$ is representation-finite. By Theorem~\ref{domestic}(1), each indecomposable vector bundle in $E^{\perp}$ is exceptional. It is also true that each indecomposable torsion sheaf $T$ in $E^{\perp}$ is exceptional. Or else 
   by Proposition~\ref{contain_non_exc}, $E^{\perp}$ is $\vec c$-invariant. So $E^{\perp}$ is a perpendicular category of a torsion exceptional sequence by Theorem~\ref{c_inv}, which is a contradiction. So Each indecomposable object in $E^{\perp}$ is exceptional. By \cite[Theorem~4.9]{hubery2016categorification}, $A$ is of representation-finite.
\end{proof}

\section{Noncrossing partitions of type $\tilde{A}$}

Let $Q_{n}$ be the following quiver of $\tilde{A}_n$ type. 
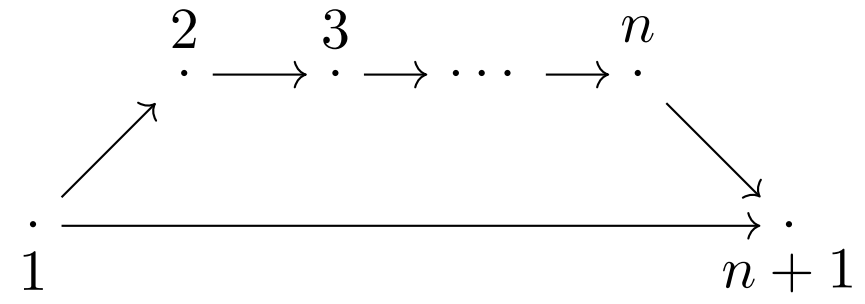
\[\tikz{
   \node (a1) at (0,0) {$\cdot$};
   \node (a2) at (1,1) {$\cdot$};
   \node (a3) at (2,1) {$\cdot$};
   \node (a4) at (3,1) {$\cdots$};
   \node (a5) at (4,1) {$\cdot$};
   \node (a6) at (5,0) {$\cdot$};

   \node (b1) at (0,-0.3) {$1$};
   \node (b2) at (1,1.3) {$2$};
   \node (b3) at (2,1.3) {$3$};
   \node (b5) at (4,1.3) {$n$};
   \node (b6) at (5,-0.3) {$n+1$};

   \graph{
           (a1)->(a2)->(a3)->(a4)->(a5)->(a6);
           (a1)->(a6);
   }
}\]
The path algebra $k Q_n$ is a hereditary algebra. Denote by
\[W_{Q_n}=\grp{s_1, s_2, \cdots, s_{n+1}}{
   \begin{aligned}
      s_i^2=1, (s_js_{j+1})^3=1, (s_1s_{n+1})^3=1, \\
      i=1,2,\cdots,n+1, j=1,2,\cdots,n
   \end{aligned}
}\]
be the Weyl group of $k Q_n$, where $s_i$ represents the simple reflection reduced by the $i$-th simple module $S_i$. Fix the standard Coxeter element 
\[
c_{Q_n}=s_1 s_2 \cdots s_{n+1},
\]
then by Theorem~\ref{hubery_krause}, there is an isomorphism of posets
\[\mathrm{Exc}(k Q_n \text{-} \mathrm{mod}) \cong \mathrm{NC}(W_{Q_n}, c_{Q_n}).\]
In fact, $\mathrm{Exc}(k Q_n \text{-} \mathrm{mod})$ is a lattice. We will prove it with the help of weighted projective lines. 

Let $\mathbb X(n)$ be a weighted projective line of weight type $\bm p=(n)$. The grading group
\[L(\bm p)=\langle \vec x_1, \vec x_2 \ | \ n \vec x_1=\vec x_2 \rangle.\]
The $L(\bm p)$-graded homogeneous coordinate ring is normalized to be $k[x_1,x_2]$ with $\mathrm{deg}(x_1)=\vec x_1$, $\mathrm{deg}(x_2)=\vec x_2$. Then the only exceptional point is $\infty$. Clearly $\mathbb X(n)$ is of domestic type. By \cite[Proposition~4.1]{geigle1987class}, there is a derived equivalence
\[\mathcal D^b(\mathrm{coh} \text{-} \mathbb X(n)) \cong \mathcal D^b(k Q_n \text{-}\mathrm{mod}).\]
Then the following proposition implies that $\mathrm{Exc}(k Q_n \text{-} \mathrm{mod})$ is a lattice, and so as $\mathrm{NC}(W_{Q_n}, c_{Q_n})$.

\begin{prop} \label{prop lattice or not}
   The following statements hold.

   (1) $\mathrm{Exc}(\mathrm{coh} \text{-} \mathbb X(n))$ is a complete lattice.

   (2) Let $\mathbb X$ be a weighted projective line with at least 2 exceptional points. Then $\mathrm{Exc}(\mathrm{coh} \text{-} \mathbb X)$ is not a lattice.
\end{prop}

\begin{proof}
   (1) We only need to show that any subset $\{\mathcal W_i \}_{i \in I}$ of $\mathrm{Exc}(\mathrm{coh} \text{-} \mathbb X(n))$ has a least upper bound. This implies that any subset $S$ of $\mathrm{Exc}(\mathrm{coh} \text{-} \mathbb X(n))$ has the greatest lower bound since the lower bounds of $S$ has a least upper bound. So $\mathrm{Exc}(\mathrm{coh} \text{-} \mathbb X(n))$ is a complete lattice. Note that lower bounds of $S$ must exist since the zero subcategory is a lower bound of $S$.
   
   If there is a vector bundle in some $\mathcal W_i$, then $\langle \mathcal W_i \rangle$ contains a vector bundle, by Proposition~\ref{thick_domestic_bundle}, $\langle \mathcal W_i \rangle_{i \in I} \in \mathrm{Exc}(\mathrm{coh} \text{-} \mathbb X(n))$. If each $\mathcal W_i$ is formed by torsion sheaves, then there are two possibilities: 
   \[
      \langle \mathcal W_i \rangle_{i \in I} \in \mathrm{Exc}(\mathrm{tor}_{\infty} \text{-} \mathbb X(n)) \text{ or } \langle \mathcal W_i \rangle_{i \in I} \in \mathrm{NExc}(\mathrm{tor}_{\infty} \text{-} \mathbb X(n)).
   \]
   If the former holds, then there is nothing to prove. If the latter holds, by Theorem~\ref{tube perpendicular}, there exists a torsion exceptional sequence $(E_1, \cdots, E_m)$ such that 
   \[
      \langle \mathcal W_i \rangle_{i \in I}=(E_1, \cdots, E_m)^{\perp} \cap \mathrm{tor}_{\infty} \text{-} \mathbb X(n).
   \]
   Then $(E_1, \cdots, E_m)^{\perp}$ is a upper bound of $\{\mathcal W_i \}_{i \in I}$. By Proposition~\ref{contain_non_exc}, any upper bound of $\{\mathcal W_i \}_{i \in I}$ is $\vec c$-invariant, and by Theorem~\ref{c_inv}, it has to be a perpendicular category of some exceptional sequence. So $(E_1, \cdots, E_m)^{\perp}$ is the least upper bound of $\{\mathcal W_i \}_{i \in I}$.

   (2) Take $\mathcal W_1=\langle S_{1,0} \rangle$, $\mathcal W_2=\langle S_{1,1},\cdots,S_{1,p_1-1} \rangle$. We will show that $\mathcal W_1$ and $\mathcal W_2$ has no least upper bounds in $\mathrm{Exc}(\mathrm{coh} \text{-} \mathbb X)$. Or else, the least upper bound $\mathcal W$ contains a vector bundle since $\langle \mathcal W_1, \mathcal W_2 \rangle=\mathrm{tor}_{\lambda_1} \text{-} \mathbb X$. So any two upper bounds of $\mathcal W_1$ and $\mathcal W_2$ have a common vector bundle. But it is not true. By Theorem~\ref{c_inv}, there exist $\vec c$-invariant wide subcategories $\mathcal V_1 \supseteq \mathcal W_1$ and $\mathcal V_2 \supseteq \mathcal W_2$ such that 
   \[
        \mathcal V_1 \cap \mathrm{tor}_{\lambda_1} \text{-} \mathbb X=\mathrm{tor}_{\lambda_1} \text{-} \mathbb X \text{\ \ and \ \ } \mathcal V_1 \cap \mathrm{tor}_{\lambda_2}\text{-} \mathbb X=\langle S_{2,1}^{[p_2]} \rangle,
   \]
   \[
        \mathcal V_2 \cap \mathrm{tor}_{\lambda_1} \text{-} \mathbb X=\mathrm{tor}_{\lambda_1} \text{-} \mathbb X \text{\ \ and \ \ } \mathcal V_2 \cap \mathrm{tor}_{\lambda_2}\text{-} \mathbb X=\langle S_{2,2}^{[p_2]} \rangle.
   \]
   So 
   \[
      \mathcal V_1 \cap \mathcal V_2 \cap \mathrm{tor}_{\lambda_2}\text{-} \mathbb X=0 \in \mathrm{Exc}(\mathrm{tor}_{\lambda_2}\text{-} \mathbb X).
   \]
   By Theorem~\ref{c_inv}, $\mathcal V_1 \cap \mathcal V_2$ only formed by torsion sheaves, or else 
   \[
      \mathcal V_1 \cap \mathcal V_2 \cap \mathrm{tor}_{\lambda_2}\text{-} \mathbb X \in \mathrm{NExc}(\mathrm{tor}_{\lambda_2}\text{-} \mathbb X).
   \]
   This is a contradiction. 
\end{proof}

\begin{rmk}
   Not all noncrossing partitions of type $\tilde A$ are lattices. For the Weyl group $W(Q_n)$, take a new Coxeter element $c'=s_1 s_2 \cdots s_{n-1} s_{n+1} s_n$. Then the poset $\mathrm{NC}(W(Q_n),c')$ is isomorphic to $\mathrm{Exc}(\mathrm{coh} \text{-} \mathbb X(n-1,2))$. By Proposition~\ref{prop lattice or not}(2), it is not a lattice. In fact, for a poset of noncrossing partitions of any affine type, \cite[Theorem 8.9]{mccammond2017artin} extends the poset into a lattice, which is the crystallographic interval in the corresponding crystallographic group. Now it is usually called McCammond-Sulway lattice. 
\end{rmk}

Let $S$ be a poset and $S'$ its subset. We denote by $\mathrm{Upp}(S')$ the set of upper bounds of $S'$, and denote by $\mathrm{Low}(S')$ the set of lower bounds of $S'$. For a weighted projective line $\mathbb X$ of domestic type, We give a poset-theoretical description of $\mathrm{Exc}(\mathrm{coh}_0 \text{-} \mathbb X)$.

\begin{prop} \label{prop Exc(tor) poset-theoretical}
   Let $\mathbb X$ be a weighted projective line of domestic type, and let $\mathcal W \in \mathrm{Exc}(\mathrm{coh} \text{-} \mathbb X)$. Then $\mathcal W \in \mathrm{Exc}(\mathrm{coh}_0 \text{-} \mathbb X)$ if and only if $\mathrm{Upp}(\mathcal W)$ is an infinite set.
\end{prop}

\begin{proof}
   If $\mathcal W \in \mathrm{Exc}(\mathrm{coh}_0 \text{-} \mathbb X)$, then by Proposition~\ref{prop tor perpendicular}, $\mathcal W^{\perp}$ has a direct summand which is equivalent to $\mathrm{coh} \text{-} \mathbb X'$, where $\mathbb X'$ is a weighted projective line of a smaller weight type. It is clear that $\mathrm{Exc}(\mathrm{coh} \text{-} \mathbb X')$ is an infinite set. Since taking right perpendicular category is a duality on the poset $\mathrm{Exc}(\mathrm{coh} \text{-} \mathbb X)$, there is infinitely $\mathcal V \in \mathrm{Exc}(\mathrm{coh} \text{-} \mathbb X)$ such that $\mathcal W \subseteq \mathcal V$, i.e. the set $\mathrm{Upp}(\mathcal W)$ is infinite.

   If $\mathcal W \in \mathrm{Exc}(\mathrm{coh} \text{-} \mathbb X)$ containing a vector bundle, we want to prove that there are only finitely wide subcategories in $\mathrm{Exc}(\mathrm{coh} \text{-} \mathbb X)$ containing $\mathcal W$. Since $\mathbb X$ is of domestic type, by Theorem~\ref{domestic}(1), $\mathcal W$ contains an exceptional vector bundle $E$. We claim that $\mathrm{Upp}(E)$ is a finite set. 
   It is equivalent to say $\mathrm{Low}(E^{\perp})$ is a finite set. By Theorem~\ref{prop domestic rep-finite}, the category $E^{\perp}$ is equivalent to the category of finite presented modules over a representation-finite algebra, so $\mathrm{Exc}(E^{\perp})$ is a finite set, i.e. $\mathrm{Low}(E^{\perp})$ is a finite set. It follows that $\mathrm{Upp}(\mathcal W) \subseteq \mathrm{Upp}(E)$ is a finite set.

   The proof is complete.
\end{proof}

We define a subset 
\[B(Q_n):=\{w \in \mathrm{NC}(W_{Q_n},c_{Q_n}) \ | \ |\mathrm{Upp}(w) \cup \mathrm{Low}(w)|=\infty \}\]
of $\mathrm{NC}(W_{Q_n},c_{Q_n})$. It is isomorphic to $\mathrm{NC}^B(n)$ and corresponds to $\mathrm{Exc}(\mathrm{coh} \text{-} \mathbb X(n))^{\vec c}$. See the next theorem.

\begin{thm} \label{thm B_n embedded in A_n tilde}
   The following statements hold.

   (1) There is a commutative diagram of lattices
   \[
   \xymatrix{\mathrm{Exc}(\mathrm{coh} \text{-} \mathbb X(n))^{\vec c} \ar@{^(->}[r] \ar[d]^{\cong} & \mathrm{Exc}(\mathrm{coh} \text{-} \mathbb X(n)) \ar[d]^{\cong} \\
   B(Q_n) \ar@{^(->}[r] & \mathrm{NC}(W_{Q_n},c_{Q_n}).
   }
   \]

   (2) There is an embedding of lattices 
   \[
     i: \mathrm{NC}^B(n) \hookrightarrow \mathrm{NC}(W_{Q_n},c_{Q_n})
   \] 
   which preserves Kreweras complement in the two lattices separately.
\end{thm}

\begin{proof}
   (1) We already have the standard isomorphism of lattices
   \[\mathrm{cox}: \mathrm{Exc}(\mathrm{coh} \text{-} \mathbb X(n)) \cong \mathrm{NC}(W_{Q_n},c_{Q_n}). \]
   So we only need to prove that 
   \[ \mathrm{cox} (\mathrm{Exc}(\mathrm{coh} \text{-} \mathbb X(n))^{\vec c}) =B(Q_n). \]
   By Theorem~\ref{c_inv}, 
   \[\mathrm{Exc}(\mathrm{coh} \text{-} \mathbb X(n))^{\vec c}=\mathrm{Exc}(\mathrm{tor}_{\infty} \text{-} \mathbb X(n)) \sqcup \mathrm{Exc}(\mathrm{tor}_{\infty} \text{-} \mathbb X(n))^{\perp},\] 
   where 
   \[\mathrm{Exc}(\mathrm{tor}_{\infty} \text{-} \mathbb X(n))^{\perp}:=\{\mathcal W^{\perp} \ | \ \mathcal W \in \mathrm{Exc}(\mathrm{tor}_{\infty} \text{-} \mathbb X(n)) \}.\]
   By Proposition~\ref{prop Exc(tor) poset-theoretical}, $\mathcal W \in \mathrm{Exc}(\mathrm{tor}_{\infty} \text{-} \mathbb X(n))$ if and only if $|\mathrm{Upp}(\mathcal W)|=\infty$. It follows that $\mathcal W \in \mathrm{Exc}(\mathrm{tor}_{\infty} \text{-} \mathbb X(n))^{\perp}$ if and only if $|\mathrm{Low}(\mathcal W)|=\infty$. So 
   \[\mathcal W \in \mathrm{Exc}(\mathrm{coh} \text{-} \mathbb X(n))^{\vec c} \text{ if and only if } |\mathrm{Upp}(\mathcal W) \cap \mathrm{Low}(\mathcal W)|=\infty.\]
   And then by the definition of $B(Q_n)$, we have $\mathrm{cox} (\mathrm{Exc}(\mathrm{coh} \text{-} \mathbb X(n))^{\vec c}) =B(Q_n)$.

   It remains to prove that $\mathrm{Exc}(\mathrm{coh} \text{-} \mathbb X(n))^{\vec c}$ is a lattice. In fact, since 
   \[\mathrm{Exc}(\mathrm{coh} \text{-} \mathbb X(n))^{\vec c}=\mathrm{Exc}(\mathrm{tor}_{\infty} \text{-} \mathbb X(n)) \sqcup \mathrm{Exc}(\mathrm{tor}_{\infty} \text{-} \mathbb X(n))^{\perp},\]
   By Theorem~\ref{c_inv}, we have an isomorphism of posets
   \begin{align*}
      \mathrm{Exc}(\mathrm{coh} \text{-} \mathbb X(n))^{\vec c} &\xrightarrow{\sim} \mathrm{Exc}(\mathrm{tor}_{\infty} \text{-} \mathbb X(n)) \\
      \mathcal W &\mapsto \mathcal W \cap \mathrm{tor}_{\infty} \text{-} \mathbb X(n).
   \end{align*}
   And then $\mathrm{Exc}(\mathrm{coh} \text{-} \mathbb X(n))^{\vec c} \cong \mathrm{Exc}(\mathrm{tor}_{\infty} \text{-} \mathbb X(n)) \cong \mathrm{NC}^B(n)$ is a lattice.

   (2) This embedding $i$ is the composition of the following maps
   \begin{align*}
      \mathrm{NC}^B(n) &\cong \mathrm{Exc}(\mathrm{tor}_{\infty} \text{-} \mathbb X(n)) \cong \mathrm{Exc}(\mathrm{coh} \text{-} \mathbb X(n))^{\vec c} \\ &\hookrightarrow \mathrm{Exc}(\mathrm{coh} \text{-} \mathbb X(n)) \cong \mathrm{NC}(W_{Q_n},c_{Q_n}).
   \end{align*}
   By Theorem~\ref{thm tube B}, the isomorphism 
   \[f: \mathrm{NC}^B(n) \cong \mathrm{Exc}(\mathrm{tor}_{\infty} \text{-} \mathbb X(n))\]
   satisfies that 
   \[f(K_B(P))=\prescript{\perp}{}{f(P)}, \ \forall P \in \mathrm{NC}^B(n),\]
   where $K_B$ denotes the Kreweras complement in the lattice $\mathrm{NC}^B(n)$. And then by the construction of $\mathrm{Exc}(\mathrm{tor}_{\infty} \text{-} \mathbb X(n)) \cong \mathrm{Exc}(\mathrm{coh} \text{-} \mathbb X(n))^{\vec c}$ in (1), the isomorphism 
   \[g: \mathrm{NC}^B(n) \cong \mathrm{Exc}(\mathrm{coh} \text{-} \mathbb X(n))^{\vec c}, \text{ satisfies } g(K_B(P))=\prescript{\perp}{}{g(P)}.\]
   That means the embedding
   \[h: \mathrm{NC}^B(n) \hookrightarrow \mathrm{Exc}(\mathrm{coh} \text{-} \mathbb X(n)), \text{ satisfying } h(K_B(P))=\prescript{\perp}{}{h(P)}.\]
   Finally, since the isomorphism 
   \[\mathrm{cox}: \mathrm{Exc}(\mathrm{coh} \text{-} \mathbb X(n)) \cong \mathrm{NC}(W_{Q_n},c_{Q_n}) \text{ satisfying } K_{\tilde A}(\mathrm{cox}(\mathcal W))=\mathrm{\prescript{\perp}{}{\mathcal W}},\]
   where $K_{\tilde A}$ denotes Kreweras complement in $\mathrm{NC}(W_{Q_n},c_{Q_n})$. So we have 
   \[
     i(K_B(P))=K_{\tilde A}(i(P)), \ \forall P \in \mathrm{NC}^B(n).
   \]
   That is, the embedding $i: \mathrm{NC}^B(n) \hookrightarrow \mathrm{NC}(W_{Q_n},c_{Q_n})$ preserves Kreweras complement.
\end{proof}

\begin{rmk}
   The existence of inclusions from noncrossing partitions of type $B$ to noncrossing partitions of any affine type can be concluded by \cite[Proposition~7.6]{mccammond2017artin}. It is well known that affine type noncrossing partitions corresponds to finite generated wide subcategories of $\mathrm{coh}\text{-} \mathbb X$ of domestic type. In fact, from our proof, it is easy to see that for a weighted projective line $\mathbb X(p_1, \cdots, p_n)$ of any type, there are inclusions 
   \[i_k: \mathrm{NC}^B(p_k) \hookrightarrow \mathrm{Exc}(\mathrm{coh} \text{-} \mathbb X), \ k=1,2, \cdots,n\] 
   which preserve ``Kreweras complement", i.e. 
   \[i_k(K_B(P))=\prescript{\perp}{}{i_k(P)}, \ k=1,2, \cdots,n.\]
\end{rmk}

We have the following observation through this categorification.

\begin{prop} \label{prop clamp}
   Let $n \geq 2$, and let $w \in \mathrm{NC}(W_{Q_n},c_{Q_n})$. Then there exists a nontrivial element (i.e. neither the greatest nor the least) $v \in B(Q_n)$ such that $w \leq v$ or $w \geq v$.
\end{prop}

\begin{proof}
   By Theorem~\ref{thm B_n embedded in A_n tilde}, it is equivalent to prove that for any wide subcategory $\mathcal W \in \mathrm{Exc}(\mathrm{coh} \text{-} \mathbb X(n))$, there exists a nontrivial wide subcategory $\mathcal V \in \mathrm{Exc}(\mathrm{coh} \text{-} \mathbb X(n))^{\vec c}$ such that $\mathcal W \subseteq \mathcal V$ or $\mathcal W \supseteq \mathcal V$. If $\mathcal W$ is $\vec c$-invariant, then there is nothing to prove. Assume that $\mathcal W$ is not $\vec c$-invariant, then $\mathcal W$ contains a vector bundle, and then by Theorem~\ref{c_inv}, $\mathcal W \cap \mathrm{tor}_{\infty} \text{-} \mathbb X(n) \in \mathrm{Exc}(\mathrm{tor}_{\infty} \text{-} \mathbb X(n))$. Now that there exists a nontrivial $\mathcal V \in \mathrm{Exc}(\mathrm{tor}_{\infty} \text{-} \mathbb X(n))$ such that $\mathcal W \supseteq \mathcal V$ if and only if $\mathcal W \cap \mathrm{tor}_{\infty} \text{-} \mathbb X(n) \neq 0$. Since taking right perpendicular category gives dualities on $\mathrm{Exc}(\mathrm{coh} \text{-} \mathbb X(n))$ and $\mathrm{Exc}(\mathrm{coh} \text{-} \mathbb X(n))^{\vec c}$, we only need to prove the following: Let $\mathcal W \in \mathrm{Exc}(\mathrm{coh} \text{-} \mathbb X(n))$ and not $\vec c$-invariant, then 
   \[\mathcal W \cap \mathrm{tor}_{\infty} \text{-} \mathbb X(n)=0 \text{ implies that } \mathcal W^{\perp} \cap \mathrm{tor}_{\infty} \text{-} \mathbb X(n) \neq 0.\]

   By Theorem~\ref{domestic}(2), line bundles $\{\mathcal O, \mathcal O(\vec{x}_1), \cdots, \mathcal{O}((n-1)\vec x_1), \mathcal O(\vec{x}_2)\}$ form a $\tau$-slice of the AR-quiver of $\mathrm{vect}(\mathbb X(n))$. So each indecomposable vector bundle is a line bundle. We now determine when $(\mathcal O, \mathcal O(\vec l))$ forms an exceptional pair. 
   \[
   \mathrm{Hom}_{\mathbb X}(\mathcal O(\vec l), \mathcal O)=0 \text{ if and only if } 0 \nleq -l, 
   \]
   if and only if $-l \leq \vec c+\vec \omega=-\vec x_1$. By Serre duality, 
   \[
      \mathrm{Ext}_{\mathbb X}^1(\mathcal O(\vec l), \mathcal O)=0 \text{ if and only if } \mathrm{Hom}_{\mathbb X}(\mathcal O, \mathcal O(\vec l+\vec \omega))=0,
   \]
   if and only if $0 \nleq \vec l+\vec \omega$, if and only if $\vec l+\vec \omega \leq \vec c+\vec \omega$. In conclusion, $(\mathcal O, \mathcal O(\vec l))$ forms an exceptional pair if and only if $\vec x_1 \leq \vec l \leq \vec c$. 

   If $\mathcal W \cap \mathrm{tor}_{\infty} \text{-} \mathbb X(n)=\mathcal W^{\perp} \cap \mathrm{tor}_{\infty} \text{-} \mathbb X(n)=0$, then all indecomposable objects in $\mathcal W$ or $\mathcal W^{\perp}$ are line bundles. Since $n \geq 2$, it follows that $\mathcal W$ or $\mathcal W^{\perp}$ contains an exceptional pair of the form $(\mathcal O(\vec x), \mathcal O(\vec x+m \vec x_1))$, where $\vec x_1 \leq \vec l \leq \vec c$. The following exact sequence 
   \[
     0 \to \mathcal O(\vec x) \to \mathcal O(\vec x+m\vec x_1) \to S_{1,m}^{[m]} \to 0
   \]
   implies that $S_{1,m}^{[m]} \in \mathcal W \text{ or }\mathcal W^{\perp}$, which contradicts to 
   \[
      \mathcal W \cap \mathrm{tor}_{\infty} \text{-} \mathbb X(n)=\mathcal W^{\perp} \cap \mathrm{tor}_{\infty} \text{-} \mathbb X(n)=0.
   \]
   The proof is complete.
\end{proof}

\begin{rmk}
	In general, let $\mathbb X$ be a weighted projective line of any domestic type, $\mathcal W$ a wide subcategory which is not $\vec c$-invariant. Then $\mathcal W \cap \mathrm{coh}_0 \text{-} \mathbb X=0$ does not imply $\mathcal W^\perp \cap \mathrm{coh}_0 \text{-} \mathbb X \neq 0$. For example, let $\mathbb X$ be of type $(2,2,2)$ and let $E$ be a rank 2 indecomposable bundle in $\mathrm{coh}\text{-} \mathbb X$. Then either $\langle E \rangle$ or $E^\perp$ only contains vector bundles.
\end{rmk}

\begin{exa}
   We describe the lattice $\mathrm{Exc}(\mathrm{coh} \text{-} \mathbb X(2))$. Denote by
   \[\mathcal T_0=\langle \mathcal O \rangle, \ \mathcal T_1=\langle \mathcal O, \mathcal O(\vec x_1) \rangle, \ \mathcal T_2=\langle \mathcal O, \mathcal O(\vec x_2)\rangle.\] 
   From the proof of Proposition~\ref{prop clamp}, it is straightforward that if $(\mathcal O, \mathcal O(\vec l))$ is an exceptional pair, then either $\vec l=\vec x_1$ or $\vec l=\vec x_2$. Notice that 
   \[
      \langle \mathcal O,S_{1,0} \rangle=\langle \mathcal O, \mathcal O(-\vec x_1) \rangle, \ \langle \mathcal O,S_{1,1} \rangle=\langle \mathcal O, \mathcal O(\vec x_1) \rangle,
   \]
   so finitely generated wide subcategories of rank 2 is of the form $\mathcal T_1(\vec l)$ or $\mathcal T_2 (\vec l)$ for some $\vec l \in L(\bm p)$. It follows that 
   \begin{align*}
      \mathrm{Exc}(\mathrm{coh} \text{-} \mathbb X)= & \{\mathcal T_i (n \vec x_1) \ | \ i=0,1,2; n \in \mathbb Z \} \\
      & \cup \{0, \mathrm{coh} \text{-} \mathbb X, \langle S_{1,0} \rangle, \langle S_{1,1} \rangle \},
   \end{align*}
   where $S_{1,i}$ is defined to be the cokernel of $\mathcal O((i-1)\vec x_1) \xrightarrow{x_1} \mathcal O(i \vec x_1)$. And the poset structure is also clear: 
   \begin{gather*}
      \langle S_{1,1} \rangle <\mathcal T_1(2n \vec x_1), \ \langle S_{1,0} \rangle <\mathcal T_1((2n-1) \vec x_1), \ \mathcal T_0(n \vec x_1)< \mathcal T_1(n \vec x_1), \\
      \mathcal T_0(n \vec x_1)< \mathcal T_1((n-1) \vec x_1), \ \mathcal T_0(2n \vec x_1)<\ \mathcal T_2 , \ \mathcal T_0((2n-1) \vec x_1)<\mathcal T_2 (\vec x_1),
   \end{gather*}
   for all $n \in \mathbb Z$. 
   
   the Hasse diagram of $\mathrm{Exc}(\mathrm{coh} \text{-} \mathbb X(2))$ is as follows. The lattice formed by black points corresponds to $\mathrm{Exc}(\mathrm{coh} \text{-} \mathbb X(2))^{\vec c}$, which is isomorphic to $\mathrm{NC}^B(2)$. As is stated in Proposition~\ref{prop Exc(tor) poset-theoretical}, there are infinitely points linking to black points, and there are only finite points linking to white points. And As is stated in Proposition~\ref{prop clamp}, each white point links to some nontrivial black points.
   
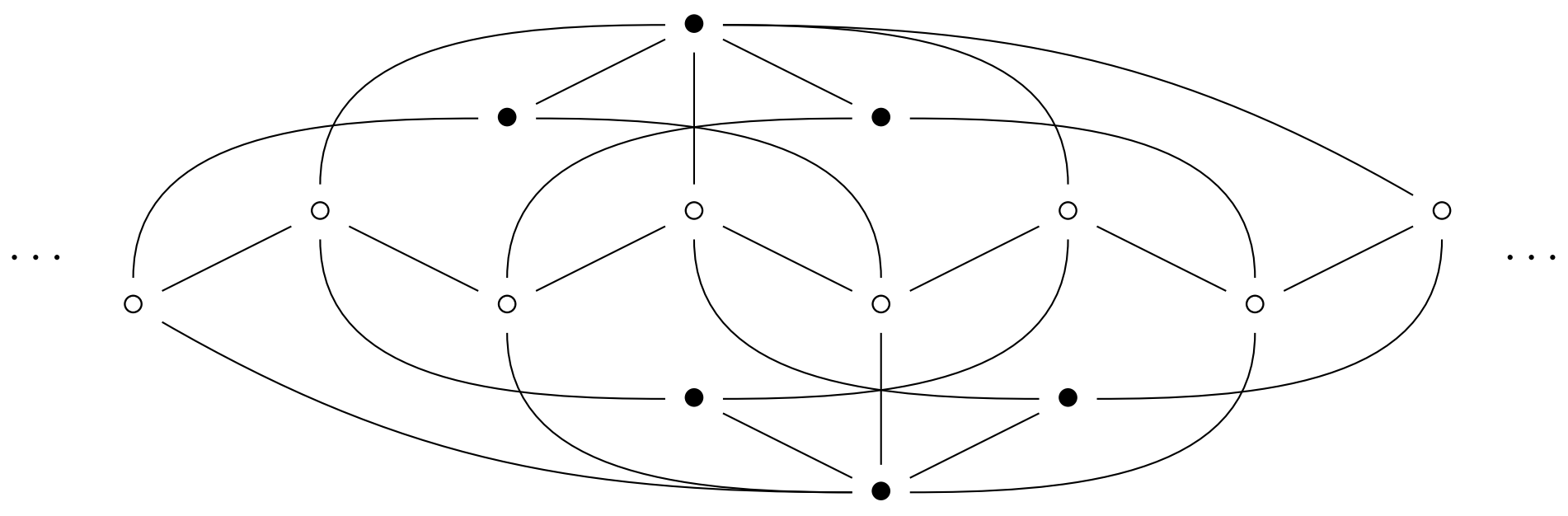
\[
   \tikz[scale=0.75]{
      \node (a1) at (0,0) {$\circ$};
      \node (a2) at (2,1) {$\circ$};
      \node (a3) at (4,0) {$\circ$};
      \node (a4) at (6,1) {$\circ$};
      \node (a5) at (8,0) {$\circ$};
      \node (a6) at (10,1) {$\circ$};
      \node (a7) at (12,0) {$\circ$};
      \node (a8) at (14,1) {$\circ$};

      \node (b1) at (6,3) {$\bullet$};
      \node (b2) at (4,2) {$\bullet$};
      \node (b3) at (8,2) {$\bullet$};

      \node (c1) at (8,-2) {$\bullet$};
      \node (c2) at (6,-1) {$\bullet$};
      \node (c3) at (10,-1) {$\bullet$};

      \node (d1) at (-1,0.5) {$\cdots$};
      \node (d2) at (15,0.5) {$\cdots$};

      \graph{
           (a1)--(a2)--(a3)--(a4)--(a5)--(a6)--(a7)--(a8);
           (a1)--[out=90, in=180](b2)--[out=0, in=90](a5);
           (a3)--[out=90, in=180](b3)--[out=0, in=90](a7);
           (a2)--[out=90, in=180](b1)--[out=0, in=90](a6);
           (a4)--(b1)--[out=0, in=150](a8);
           (b2)--(b1)--(b3);
           (a2)--[out=-90, in=180](c2)--[out=0, in=-90](a6);
           (a4)--[out=-90, in=180](c3)--[out=0, in=-90](a8);
           (a3)--[out=-90, in=180](c1)--[out=0, in=-90](a7);
           (a1)--[out=-30, in=180](c1)--(a5);
           (c2)--(c1)--(c3);
         }
   }
   \]
\end{exa}
For the convenience of the reader, we present the corresponding wide subcategories in the Hasse diagram as follows.
\[
   \tikz[scale=0.75]{
      \node (a1) at (0,0) {$\mathcal T_0$};
      \node (a2) at (2,1) {$\mathcal T_1$};
      \node (a3) at (4,0) {$\mathcal T_0(\vec x_1)$};
      \node (a4) at (6,1) {$\mathcal T_1(\vec x_1)$};
      \node (a5) at (8,0) {$\mathcal T_0(2 \vec x_1)$};
      \node (a6) at (10,1) {$\mathcal T_1(\vec x_2)$};
      \node (a7) at (12,0) {$\mathcal T_0(3 \vec x_1)$};
      \node (a8) at (14,1) {$\mathcal T_1(\vec x_3)$};

      \node (b1) at (6,3) {$\mathrm{coh} \text{-} \mathbb X(2)$};
      \node (b2) at (4,2) {$\mathcal T_2$};
      \node (b3) at (8,2) {$\mathcal T_2(\vec x_1)$};

      \node (c1) at (8,-2) {$\{0\}$};
      \node (c2) at (6,-1) {$\langle S_{1,1} \rangle$};
      \node (c3) at (10,-1) {$\langle S_{1,0} \rangle$};

      \node (d1) at (-1,0.5) {$\cdots$};
      \node (d2) at (15,0.5) {$\cdots$};
   }
   \]

\section*{Acknowledgements}

The author thanks Professor Xiao-Wu Chen and Professor Shiquan Ruan for many helpful discussions and advices. The author thanks Professor Nathan Reading for constructive comments.

\end{document}